\documentclass[11pt]{article}
\usepackage{hyperref}
\usepackage{marvosym}
\usepackage{graphicx}
\usepackage{latexsym,amsmath,amsfonts,amscd, amsthm, dsfont}
\usepackage{bm,color}
\usepackage{epsfig,verbatim,epstopdf,graphics}
\usepackage{subfigure}
\usepackage{changebar}
\usepackage{multirow}

\usepackage{algorithmic}

\usepackage{yhmath}
 \usepackage{booktabs} 
\usepackage{tikz}
\usepackage{verbatim}
\usepackage{diagbox}
\usetikzlibrary{arrows,backgrounds,snakes,shapes}
 \numberwithin{equation}{section}

\graphicspath{{./}{./figure/}}
\allowdisplaybreaks

\newtheoremstyle{plainNoItalics}{}{}{\normalfont}{}{\bfseries}{.}{ }{}

\theoremstyle{plain}
\newtheorem{thm}{Theorem}[section]

\theoremstyle{plainNoItalics}

\newtheorem{defn}[thm]{Definition}
\newtheorem{rem}[thm]{Remark}

\newtheorem{exa}[thm]{Example}

\newcommand{\bit}{\begin{itemize}}
\newcommand{\eit}{\end{itemize}}
\newcommand{\beq}{\begin{equation}}
\newcommand{\eeq}{\end{equation}}
\newcommand{\be}{\begin{eqnarray}}
\newcommand{\ee}{\end{eqnarray}}
\newcommand{\beno}{\begin{eqnarray*}}
\newcommand{\eeno}{\end{eqnarray*}}

\newcommand{\email}[1]{\protect\href{mailto:#1}{#1}}

\makeatletter

\newcommand{\Rmnum}[1]{\expandafter\@slowromancap\romannumeral #1@}

\makeatother

\begin{document}

\baselineskip=1.8pc

\title{A mass conservative Eulerian-Lagrangian Runge-Kutta discontinuous Galerkin method for wave equations with large time stepping.%
  \thanks{
  Research of the first author is supported by the China Scholarship Council for 2 years' study at the University of Delaware, CSC NO. 201906340216.
Research of the second author is supported by NSF grant NSF-DMS-2111253, Air Force Office of Scientific Research FA9550-18-1-0257 and University of Delaware.
 }}
\author{Xue Hong
  \thanks{School of Mathematical Sciences, University of Science and Technology of China, Hefei, Anhui, 230026, P.R. China. (\email{xuehong1@mail.ustc.edu.cn}).}
  \and
  Jing-Mei Qiu%
  \thanks{Department of Mathematical Sciences, University of Delaware, Newark, DE, 19716, USA. (\email{jingqiu@udel.edu}).}
  }

\maketitle

 \begin{abstract}
We propose an Eulerian-Lagrangian (EL) Runge-Kutta (RK) discontinuous Galerkin (DG) method for wave equations. The method is designed based on the ELDG method for transport problems [J. Comput. Phy. 446: 110632, 2021.], which tracks solution along approximations to characteristics in the DG framework, allowing extra large time stepping sizes with stability. The wave equation can be written as a first order hyperbolic system. Considering each characteristic family, a straightforward application of ELDG will be to transform to the characteristic variables, evolve them on associated characteristic related space-time regions, and transform them back to the original variables. However, the mass conservation could not be guaranteed in a general setting. In this paper, we formulate a mass conservative semi-discrete ELDG method by decomposing each variable into two parts, each of them associated with a different characteristic family. As a result, four different quantities are evolved in EL fashion and recombined to update the solution. The fully discrete scheme is formulated by using method-of-lines RK methods, with intermediate RK solutions updated on the background mesh. Numerical results on 1D and 2D wave equations are presented to demonstrate the performance of the proposed ELDG method. These include the high order spatial and temporal accuracy, stability with extra large time stepping size, and mass conservative property.

 \end{abstract}


%

\section{Introduction}
In this paper, we propose an Eulerian-Lagrangian (EL) discontinuous Galerkin (DG) method for the wave equation in the form of
\begin{equation}
u_{tt} - \nabla \cdot( \mathbf{Q}( \mathbf{x},t ) \nabla u ) = 0, \ (\mathbf{x},t)\in\mathbb{R}^d\times[0,T],
\label{general_linearwave}
\end{equation}
where $d$ is the spatial dimension, $u: \mathbb{R}^d  \times[0,T]\rightarrow  \mathbb{R}$, and $\mathbf{Q}( \mathbf{x},t  ) = diag( q_1(\mathbf{x},t ), \cdots,q_d(\mathbf{x},t ) )$. 
Such a model could come from many fields of engineering, e.g.
petroleum engineering, defense industry, telecommunication and geoscience (\cite{durran2013numerical}, \cite{kampanis2008effective}). There are many efficient DG methods for solving the second order wave equation.
These include the interior penalty DG (IPDG) method\cite{riviere2003discontinuous,grote2006discontinuous,grote2009optimal}, local DG method\cite{baccouch2012local,xing2013energy,chou2014optimal}.
One can also reformulate the problem as a first-order hyperbolic system
 \begin{equation}
U_t+\sum_{j=1}^d (A_j(\mathbf{x},t)U)_{x_j} = 0, \ (\mathbf{x},t)\in\mathbb{R}^d\times[0,T],
\label{general_hypersystem}
\end{equation}
for which DG \cite{cockburn1989tvb} and space-time DG \cite{falk1999explicit,monk2005discontinuous} methods are available.

One effective numerical approach for solving the hyperbolic equation is the characteristic method, which updates time-dependent solution by tracking characteristics. For scalar hyperbolic equations, there have been many pioneering works.
In \cite{celia1990eulerian}, Celia etc. developed an Eulerian Lagrangian Localized Adjoint Methods (ELLAM) \cite{celia1990eulerian}, which introduces an adjoint problem for each test function in the continuous finite element framework and has been applied to different problems \cite{wang1999family,russell2002overview}. On the other hand, the ELDG \cite{wang2007eulerian,cai2020eldg} and SLDG \cite{cai2016high} are being developed in the discontinuous Galerkin finite element framework with a similar introduction of adjoint problems for test functions. Another line of development, that is closely related to this work, is the Arbitrary Lagrangian Eulerian (ALE) DG method \cite{klingenberg2017arbitrary,hong2020ALEDG_singular}. Both ELDG and ALE-DG evolve the DG solution on a dynamic moving mesh: the dynamic mesh movement of the ELDG approximates characteristics for the potential of using larger time stepping sizes with stability, whereas the mesh movement of ALE-DG is designed for better shock resolution.

In this paper, we propose a mass conservative EL RK DG method for hyperbolic systems \eqref{general_hypersystem} with large time stepping sizes. We start from 1D cases, for which we consider characteristic variables and the associated characteristic space-time regions.
For each characteristic family, a straightforward application of ELDG will be to transform to characteristic variables, evolve them on associated space-time regions, and transform them back to original variables. However, the mass conservation could not be guaranteed in a general setting.
We decompose each variable into two parts, each of them associated with different characteristic families; as a result, four different quantities are evolved in EL fashion and recombined to update the solution. The fully discrete scheme is formulated by using method-of-lines RK methods, with intermediate RK solutions updated on the background mesh. For 2D hyperbolic systems, characteristic variables are no longer constant along the bicharacteristics, characteristic Galerkin \cite{ostkamp1997multidimensional} or evolution Galerkin \cite{butler1960numerical,lax1959systems,lukavcova2000evolution} methods have been proposed by taking into account information propagated in all bicharacteristic directions. However, the algorithm implementation is very complex. In this paper, we use the splitting method for higher dimensional problem, which maintain the simplicity of ELDG method for 1D cases, as well as other great properties such as mass conservation, high order spatial and temporal accuracy, and allows for extra large time steps with stability.

This paper is organized as follows. In Section \ref{section:1d}, we review the ELDG for one-dimensional (1D) linear transport problems. In Section \ref{section:2d}, we develop the ELDG method for a first-order hyperbolic system by evolving each component associated with each characteristic families, recombining them to update the solution. We also develop it to 2D problems by dimensional splitting. Mass conservation of ELDG schemes are proved.
In Section \ref{section:numerical}, performance of the proposed ELDG method is shown through extensive
numerical tests. Finally, concluding remarks are made in Section \ref{section:conclusion}.

\section{Review of ELDG formulation for 1D linear transport problems} \label{section:1d}

To illustrate the key idea of the ELDG scheme, we start from a 1D linear transport equation in the following form
\begin{equation}
u_t+(a(x, t)u)_x = 0, \quad x\in[x_a, x_b].
\label{scalar1d}
\end{equation}
For simplicity, we assume periodic boundary conditions, and the velocity field $a(x, t)$ is a continuous function of space and time. We perform a partition of the computational domain $x_a=x_{\frac12}< x_{\frac32}<\cdots< x_{N+\frac12} =x_b$ as the background mesh partition.
Let $I_j=[x_{j-\frac12}, x_{j+\frac12} ]$ denote an element of length $\Delta x_j=x_{j+\frac12}-x_{j-\frac12}$ and define $\Delta x=\max_{j}\Delta x_j.$
We define the finite dimensional approximation space, $V_h^k = \{ v_h:  v_h|_{I_j} \in P^k(I_j) \}$, where $P^k(I_j)$ denotes the set of polynomials of degree at most $k$. We let $t^n$ be the $n$-th time level and $\Delta t =t^{n+1}-t^n$ to be the time-stepping size.

\begin{figure}[h!]
\centering
\begin{tikzpicture}[x=1cm,y=1cm]
  \begin{scope}[thick]

  \draw[fill=blue!5] (0.,2) -- (-1.5,0) -- (1.6,0) -- (2.4,2)
      -- cycle;
   \node[blue!70, rotate=0] (a) at ( 1. ,1.5) {\LARGE $\Omega_j$ };

   \draw (-3,3) node[fill=white] {};
    \draw (-3,-1) node[fill=white] {};
    \draw[black]                   (-3,0) node[left] {} -- (3,0)
                                        node[right]{$t^{n}$};
    \draw[black] (-3,2) node[left] {$$} -- (3,2)
                                        node[right]{$t^{n+1}$};

     \draw[snake=ticks,segment length=2.4cm] (-2.4,2) -- (0,2) node[above left] {$x_{j-\frac12}$};
     \draw[snake=ticks,segment length=2.4cm] (0,2) -- (2.4,2) node[above right] {$x_{j+\frac12}$};

          \draw[snake=ticks,segment length=2.4cm] (-2.4,0) -- (0,0);
     \draw[snake=ticks,segment length=2.4cm] (0,0) -- (2.4,0);

            \draw[blue,thick] (0.,2) node[left] {$$} -- (-1.5,0)
                                        node[black,below]{$x_{j-\frac12}^\star$ };

                \draw[blue,thick] (2.4,2) node[left] {$$} -- (1.6,0)
                                        node[black,below right]{ $x_{j+\frac12}^\star$};
\draw [decorate,color=red,decoration={brace,mirror,amplitude=9pt},xshift=0pt,yshift=0pt]
(-1.5,0) -- (1.6,0) node [red,midway,xshift=0cm,yshift=-20pt]
{\footnotesize $\widetilde{I}_j(t^{n}) = I_{j}^\star$};

\draw [decorate,color=red,decoration={brace,amplitude=10pt},xshift=0pt,yshift=0pt]
(0,2) -- (2.4,2) node [red,midway,xshift=0cm,yshift=20pt]
{\footnotesize $\widetilde{I}_j(t^{n+1}) = I_{j}$};

\draw[-latex,thick](0,1)node[right,scale=1.]{$\widetilde{I}_j(t)$}
        to[out=180,in=0] (-0.75,1) node[above left=2pt] {$$};  

\draw[-latex,thick](1,1)node[right,scale=1.]{ }
        to[out=0,in=180] (2,1) node[above left=2pt] {$$};  

  \node[blue, rotate=54] (a) at (-1.25,1.) { $\nu_{j-\frac12}$ };
  \node[blue, rotate=70] (a) at (2.35,1.) { $\nu_{j+\frac12}$ };

  \end{scope}

   \draw[-latex,dashed](1.6,1)node[right,scale=1.0]{ }
        to[out=30,in=150] (6,1) node[] {};

     \draw[-latex,red,dashed](1.6,2)node[right,scale=1.0]{ }
        to[out=40,in=140] (6,1.2) node[] {};

     \draw[-latex,red,dashed](1.3,0)node[right,scale=1.0]{ }
        to[out=-60,in=220] (6,1-0.2) node[] {};

\draw[|-|,black,thick] (6.5,0.8) node[below]{\footnotesize $\xi=x_{j-1/2}$} to (9.5,0.8)node[below]{\footnotesize $\xi=x_{j+1/2}$};
\end{tikzpicture}
\caption{Illustration for dynamic element $\widetilde{I}_j(t)$ of ELDG.}
\label{iso1}
\end{figure}
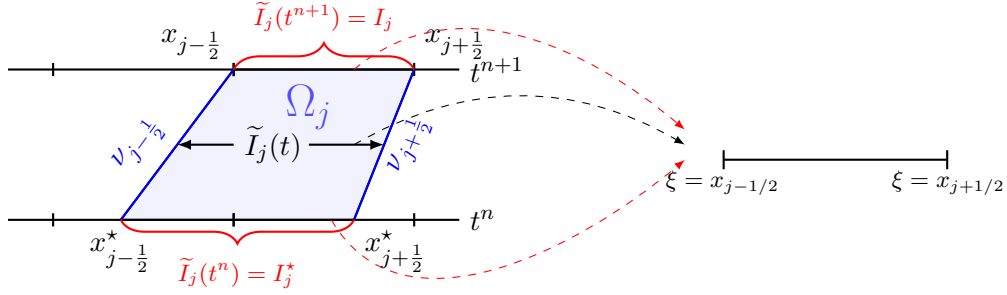

The key idea in the ELDG formulation is, design adjoint problems for test functions to take advantage of information propagation along characteristics.
The ELDG method proposed in \cite{cai2020eldg} is formulated on a space-time region $\Omega_j = \tilde{I}_j(t) \times [t^n, t^{n+1}]$ with
$$\tilde{I}_j(t) = [\tilde{x}_{j-\frac12}(t), \tilde{x}_{j+\frac12}(t)], \quad t \in [t^n, t^{n+1}]$$
being the dynamic interval, see Figure ~\ref{iso1}. Here $\tilde{x}_{j\pm\frac12}(t) = x_{j\pm\frac12} + (t-t^{n+1}) \nu_{j\pm\frac12}$ are straight lines emanating from cell boundaries $x_{j\pm\frac12}$ with slopes $\nu_{j\pm\frac12}= a(x_{j\pm\frac12},t^{n+1})$ and ${I}^\star_j \doteq \tilde{I}_j(t^n) = [x^*_{j-\frac12}, x^*_{j+\frac12}]$ is the upstream cell of $I_j$ at $t^n$. The dynamic interval of $\tilde{I}_j(t)$ can always be linearly mapped to a reference cell $\xi$ in $I_j$ by a mapping.
A local adjoint problem of \eqref{scalar1d} for the test function is defined as:
\begin{equation}
\begin{cases}
\psi_t + \alpha(x,t) \psi_x =0 ,\ (x, t)\in\Omega_j,\\
\psi(t=t^{n+1}) = \Psi(x). \quad \forall \Psi(x) \in P^k(I_j).
\end{cases}
\label{ELDGadjoint1d}
\end{equation}
Here $\alpha(x,t)$  is a bilinear function of $(x, t)$ with $\forall t\in [t^n,t^{n+1}]$,
\beq
\label{eq: a}
\alpha(x_{j-\frac12},t^{n+1}) =a(x_{j-\frac12}, t^{n+1})\doteq \nu_{j-\frac12},  \quad {\alpha}(x_{j+\frac12},t^{n+1})=a(x_{j+\frac12}, t^{n+1})\doteq \nu_{j+\frac12},
\eeq
and
\beq
\alpha(x, t)=-\nu_{j-\frac12} \frac{x-\tilde{x}_{j+\frac12}(t)}{\Delta x_j(t)} + \nu_{j+\frac12} \frac{x-\tilde{x}_{j-\frac12}(t)}{\Delta x_j(t)}\in P^1([\tilde{x}_{j-\frac12}(t), \tilde{x}_{j+\frac12}(t)]).
\label{eq: alpha_def1}
\eeq
Notice that the test function $\psi$ stays the same polynomial, if $\tilde{I}_j(t)$ is mapped to a reference interval $I_j$. The ELDG \cite{cai2020eldg} scheme can be formulated by $\int_{\tilde{I}_j(t)} \eqref{ELDGadjoint1d} \cdot u + \eqref{scalar1d} \cdot \psi $
\begin{equation}
\frac{d}{dt} \int_{\tilde{I}_j(t)}(u\psi)dx =-  \left(\hat{F} \psi\right) \left|_{\tilde{x}_{j+\frac12}(t) } \right. +   \left(\hat{F} \psi\right) \left|_{\tilde{x}_{j-\frac12}(t) } \right.   + \int_{\tilde{I}_j(t)}F\psi_xdx.
\label{mol_2}
\end{equation}
where $F(u) \doteq (a-\alpha)  u$ and $\hat{F}$ is the Lax-Friedrichs flux. A method-of-lines RK discretization can be used for high order temporal accuracy
\cite{cai2020eldg}.

\section{The ELDG algorithm for wave equation.} \label{section:2d}
The design of the ELDG algorithm for wave equation shares a similar spirit as the 1D scalar case. We start from a 1D wave equation. We first rewrite it as a linear system and formulate the ELDG scheme by tracking information along different characteristics families.
\subsection{1D wave equation}
We consider the wave equation
\begin{equation}
 u_{tt}  =  (a^2(x)u_x)_x+f(x,t).
 \label{2d_hcl}
\end{equation}
For simplicity, we assume periodic boundary conditions, and the velocity field $a(x)$ is a continuous and periodic function of space. Let $u^1=u_t$ and $u^2=u_x$, then we can rewrite ~\eqref{2d_hcl} as a linear system
\beq
\label{eq: 1d_system}
U_t+(AU)_x=F,
\eeq
where
\[
U=\begin{bmatrix}
  u^1\\
  u^2
\end{bmatrix},
\qquad
A(x)=\begin{bmatrix}
0&-a^2(x)\\
-1&0
\end{bmatrix},
\qquad
F=\begin{bmatrix}
  f(x,t)\\
  0
\end{bmatrix}.
\]
We use the following notations for the eigen-decomposition of $A(x)$:
\begin{itemize}
\item
eigenvalue:
$
\lambda^{(1)}(x), \quad \lambda^{(2)}(x),
$
where $\lambda^{(1)}(x)=a(x), \quad \lambda^{(2)}(x)=-a(x)$ for wave equation.
\item $A(x)=R(x) \Lambda R^{-1}(x)$, where $\Lambda(x)=diag(\lambda^{(1)}(x), \lambda^{(2)}(x))$,
\begin{equation}
  R(x)=\begin{bmatrix}
r_{11}(x)&r_{12}(x)\\
r_{21}(x)&r_{22}(x)
\end{bmatrix}=\begin{bmatrix}
-a(x)&a(x)\\
1&1
\end{bmatrix}
\end{equation} contains the right eigenvectors, and
\begin{equation}\label{eigenvecterandarepresent}
  R^{-1}(x)=\begin{bmatrix}
l_{11}(x)&l_{12}(x)\\
l_{21}(x)&l_{22}(x)
\end{bmatrix}=\begin{bmatrix}
\frac{-1}{2a(x)}&\frac12\\
\frac{1}{2a(x)}&\frac12
\end{bmatrix} \doteq \begin{bmatrix}
  l^{(1)}(x)\\
  l^{(2)}(x)
\end{bmatrix}
\end{equation}
contains the left eigenvectors.
\end{itemize}

In the following, we propose a mass conservative ELDG scheme for the system \eqref{eq: 1d_system} by the procedure below. The notion of its background meshes is the same as the 1D scalar case.

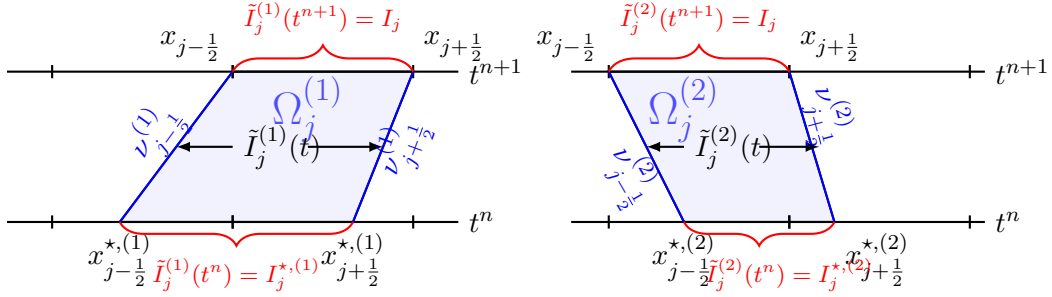
\begin{figure}[h!]
\centering
\begin{tikzpicture}[x=1cm,y=1cm]
  \begin{scope}[thick]

  \draw[fill=blue!5] (0.,2) -- (-1.5,0) -- (1.6,0) -- (2.4,2)
      -- cycle;
   \node[blue!70, rotate=0] (a) at ( 1. ,1.5) {\LARGE $\Omega^{(1)}_j$ };

   \draw (-3,3) node[fill=white] {};
    \draw (-3,-1) node[fill=white] {};
    \draw[black]                   (-3,0) node[left] {} -- (3,0)
                                        node[right]{$t^{n}$};
    \draw[black] (-3,2) node[left] {$$} -- (3,2)
                                        node[right]{$t^{n+1}$};

     \draw[snake=ticks,segment length=2.4cm] (-2.4,2) -- (0,2) node[above left] {$x_{j-\frac12}$};
     \draw[snake=ticks,segment length=2.4cm] (0,2) -- (2.4,2) node[above right] {$x_{j+\frac12}$};

          \draw[snake=ticks,segment length=2.4cm] (-2.4,0) -- (0,0);
     \draw[snake=ticks,segment length=2.4cm] (0,0) -- (2.4,0);

            \draw[blue,thick] (0.,2) node[left] {$$} -- (-1.5,0)
                                        node[black,below]{$x^{\star, (1)}_{j-\frac12}$ };

                \draw[blue,thick] (2.4,2) node[left] {$$} -- (1.6,0)
                                        node[black,below ]{ $x^{\star, (1)}_{j+\frac12}$};
\draw [decorate,color=red,decoration={brace,mirror,amplitude=9pt},xshift=0pt,yshift=0pt]
(-1.5,0) -- (1.6,0) node [red,midway,xshift=0cm,yshift=-20pt]
{\footnotesize $\tilde{I}^{(1)}_j(t^{n}) = I^{\star, (1)}_{j}$};

\draw [decorate,color=red,decoration={brace,amplitude=10pt},xshift=0pt,yshift=0pt]
(0,2) -- (2.4,2) node [red,midway,xshift=0cm,yshift=20pt]
{\footnotesize $\tilde{I}^{(1)}_j(t^{n+1}) = I_{j}$};

\draw[-latex,thick](0,1)node[right,scale=1.]{$\tilde{I}^{(1)}_j(t)$}
        to[out=180,in=0] (-0.75,1) node[above left=2pt] {$$};  

\draw[-latex,thick](1,1)node[right,scale=1.]{ }
        to[out=0,in=180] (2,1) node[above left=2pt] {$$};  

  \node[blue, rotate=54] (a) at (-1.0,1.2) { $\nu^{(1)}_{j-\frac12}$ };
  \node[blue, rotate=70] (a) at (2.25,1.) { $\nu^{(1)}_{j+\frac12}$ };




\draw[fill=blue!5] (5.,2) -- (6.,0) -- (8.,0) -- (7.4,2)
      -- cycle;
   \node[blue!70, rotate=0] (a) at ( 6. ,1.5) {\LARGE $\Omega^{(2)}_j$ };

   \draw (4.5,3) node[fill=white] {};
    \draw (4.5,-1) node[fill=white] {};
    \draw[black]                   (4.5,0) node[left] {} -- (10.,0)
                                        node[right]{$t^{n}$};
    \draw[black] (4.5,2) node[left] {$$} -- (10.,2)
                                        node[right]{$t^{n+1}$};

     \draw[snake=ticks,segment length=2.4cm] (7.4,2) -- (5.0,2) node[above left] {$x_{j-\frac12}$};
     \draw[snake=ticks,segment length=2.4cm] (9.8,2) -- (7.4,2) node[above right] {$x_{j+\frac12}$};

          \draw[snake=ticks,segment length=2.4cm] (5.,0) -- (7.4,0);
     \draw[snake=ticks,segment length=2.4cm] (7.4,0) -- (9.8,0);

            \draw[blue,thick] (5.,2) node[left] {$$} -- (6.,0)
                                        node[black,below ]{$x^{\star, (2)}_{j-\frac12}$ };

                \draw[blue,thick] (7.4,2) node[left] {$$} -- (8.,0)
                                        node[black,below right]{ $x^{\star, (2)}_{j+\frac12}$};
\draw [decorate,color=red,decoration={brace,mirror,amplitude=9pt},xshift=0pt,yshift=0pt]
(6.,0) -- (8.,0) node [red,midway,xshift=0.4cm,yshift=-20pt]
{\footnotesize $\tilde{I}^{(2)}_j(t^{n}) = I^{\star, (2)}_{j}$};

\draw [decorate,color=red,decoration={brace,amplitude=10pt},xshift=0pt,yshift=0pt]
(5,2) -- (7.4,2) node [red,midway,xshift=0cm,yshift=20pt]
{\footnotesize $\tilde{I}^{(2)}_j(t^{n+1}) = I_{j}$};

\draw[-latex,thick](6,1)node[right,scale=1.]{$\tilde{I}^{(2)}_j(t)$}
        to[out=180,in=0] (5.5,1) node[above left=2pt] {$$};  

\draw[-latex,thick](7,1)node[right,scale=1.]{ }
        to[out=0,in=180] (7.8,1) node[above left=2pt] {$$};  

  \node[blue, rotate=-60] (a) at (5.3,0.6) { $\nu^{(2)}_{j-\frac12}$ };
  \node[blue, rotate=-75] (a) at (7.9,1.4) { $\nu^{(2)}_{j+\frac12}$ };



  \end{scope}

%
%
%
\end{tikzpicture}
\caption{Illustration for dynamic elements $\widetilde{I}^{(1)}_j(t)$ (left) and $\widetilde{I}^{(2)}_j(t)$ (right) of ELDG for the first and second characteristic families of the system. }
\label{iso2}
\end{figure}

 \noindent
 {\bf (1) Two partitions of space-time regions $\Omega^{(1)}_j$ and $\Omega^{(2)}_j$.}
According to the first and second characteristic families, we partition the computational domain as two sets of space-time regions ${\Omega^{(1)}_j}$ and ${\Omega^{(2)}_j}$ respectively. Here $\Omega^{(1)}_j = \tilde{I}^{(1)}_j(t) \times [t^n, t^{n+1}]$ is related to the first characteristic family. $\tilde{I}^{(1)}_j(t) = [\tilde{x}^{(1)}_{j-\frac12}(t), \tilde{x}^{(1)}_{j+\frac12}(t)]$ is the dynamic interval emanating from cell boundaries $x_{j\pm\frac12}$ with slopes ${\nu}^{(1)}_{j\pm\frac12}$ approximating the first characteristic velocity, see Figure~\ref{iso2} (left). In general, we choose $\nu^{(1)}_{j\pm\frac12}=\lambda^{(1)}(x_{j\pm\frac12})$.
$I^{\star,(1)}_j \doteq \tilde{I}^{(1)}_j(t^n)$ is the upstream cell of $I_j$ from the first characteristic family at $t^n$.
The dynamic interval $\tilde{I}^{(1)}_j(t)$ can be linearly mapped to a reference cell $\xi\in I_j$, (see Figure~\ref{iso1}). Here, we let $\tilde{x}^{(1)}(\tau;(\xi,t^{n+1}))$ be the linear map from $\tilde{I}^{(1)}_j(t)$ to $I_j$.
Similar definition can be made to $\Omega^{(2)}_j$, $\tilde{I}^{(2)}_j(t)$ and $I^{\star, (2)}_j$
for the second characteristic family.
See Figure~\ref{iso2} (right) for illustration of $\Omega^{(2)}_j$.

\noindent {\bf (2) Adjoint Problems.}
We consider an adjoint problem for the first characteristic family on $\Omega^{(1)}_j$:
\begin{equation}
\begin{cases}
(\psi^{(1)})_t + \alpha^{(1)} (\psi^{(1)})_x =0 ,\ t\in[t^n,t^{n+1} ],\\
(\psi^{(1)})(t=t^{n+1}) = \Psi^{(1)}(x),
\end{cases}
\label{new ELDG adjointsystem1}
\end{equation}
where
\beq
\alpha^{(1)}(x, t)=-\nu^{(1)}_{j-\frac12} \frac{x-\tilde{x}^{(1)}_{j+\frac12}(t)}{\Delta x_j^{(1)}(t)} + \nu^{(1)}_{j+\frac12} \frac{x-\tilde{x}^{(1)}_{j-\frac12}(t)}{\Delta x_j^{(1)}(t)}\in P^1(\tilde{I}^{(1)}_j(t)).
\label{eq: ELDG alpha for system_def1}
\eeq
Similarly on $\Omega^{(2)}_j$:
\begin{equation}
\begin{cases}
(\psi^{(2)})_t + \alpha^{(2)} (\psi^{(2)})_x =0 ,\ t\in[t^n,t^{n+1} ],\\
(\psi^{(2)})(t=t^{n+1}) = \Psi^{(2)}(x),
\end{cases}
\label{new ELDG adjointsystem2}
\end{equation}
where
\beq
\alpha^{(2)}(x, t)=-\nu^{(2)}_{j-\frac12} \frac{x-\tilde{x}^{(2)}_{j+\frac12}(t)}{\Delta x_j^{(2)}(t)} + \nu^{(2)}_{j+\frac12} \frac{x-\tilde{x}^{(2)}_{j-\frac12}(t)}{\Delta x_j^{(2)}(t)}\in P^1(\tilde{I}^{(2)}_j(t)).
\label{eq: ELDG alpha for system_def2}
\eeq

\noindent {\bf (3) Formulation of a mass-conservative semi-discrete ELDG scheme.}
Before we present a mass-conservative semi-discrete ELDG scheme,
we formulate the scheme by a localized characteristic field. In particular, a piecewise constant $a_j$ approximating $a(x)$ in \eqref{eigenvecterandarepresent} is defined on $I_j$, and the corresponding
\begin{equation}\label{nmcapproxieigenvecterandarepresent}
  R_j=\begin{bmatrix}
r_j^{11}&r_j^{12}\\
r_j^{21}&r_j^{22}
\end{bmatrix}=\begin{bmatrix}
-a_j&a_j\\
1&1
\end{bmatrix}
\end{equation} and
\begin{equation}\label{nmcapproxieigenvectora}
  R_j^{-1}=\begin{bmatrix}
\frac{-1}{2a_j}&\frac12\\
\frac{1}{2a_j}&\frac12
\end{bmatrix} \doteq \begin{bmatrix}
  l_j^{(1)}\\
  l_j^{(2)}
\end{bmatrix}.
\end{equation}
Define $l_j^{(1)}$ is locally defined on $\Omega_j^{(1)}$ approximating $l^{(1)}(x)$.
For simplicity, we only present the first order ELDG scheme.
Take the vector product of $l^{(1)}_j$ from left with \eqref{eq: 1d_system}, we have a scalar equation
\begin{equation}\label{eqn1: equation1b}
  l^{(1)}_j( U_t+(A(x)U)_x) =l^{(1)}_j F(x,t).
\end{equation}
Next, integrating over the space-time interval $\tilde{I}^{(1)}_j(t)$, then we have
\begin{equation}
\label{nonmassconserveeqn1 of v1}
\begin{aligned}
&\frac{d}{dt}\int_{\tilde{I}^{(1)}_j(t)} (l^{(1)}_j U)dx+ l^{(1)}_j(A(x)U-\nu^{(1)}_{j+\frac{1}{2}}U)|_{\tilde{x}^{(1)}_{j+\frac12}(t)} -l^{(1)}_j(A(x)U-\nu^{(1)}_{j-\frac{1}{2}}U)|_{\tilde{x}^{(1)}_{j-\frac12}(t)}\\
&=\int_{\tilde{I}^{(1)}_j(t)} l^{(1)}_j F(x,t) dx.\\
\end{aligned}
\end{equation}
The first order ELDG discretization of eq.~\eqref{nonmassconserveeqn1 of v1} is to find $l^{(1)}_j U_h(x,t)\in P^0(\tilde{I}^{(1)}_j(t))$, so that
\begin{equation}\label{nonmassconserveeqn of v1}
\begin{aligned}
&\frac{d}{dt}\int_{\tilde{I}^{(1)}_j(t)} l^{(1)}_j U_h dx=- \left[l^{(1)}_j \widehat{(A(x)U_h-\nu^{(1)}_{j+\frac{1}{2}}U_h)}|_{j+\frac{1}{2}}\right]+ \left[l^{(1)}_{j}\widehat{(A(x)U_h-\nu^{(1)}_{j-\frac{1}{2}}U_h)}|_{j-\frac{1}{2}}\right]\\
  &+ \int_{\tilde{I}^{(1)}_j(t)} l^{(1)}_j F(x,t) dx\doteq L_{1}(U_h(t),t,\tilde{I}^{(1)}_j(t)).
\end{aligned}
\end{equation}
Here $\widehat{f^{11}}$ at a cell boundary can be taken as a monotone flux, e.g. the Lax-Friedrichs flux
\[
\widehat{(A U-\nu^{(1)}_{j+\frac{1}{2}}U)}_{{j+\frac{1}{2}}} = \frac12\left(A(x^+_{j+\frac{1}{2}})U^+_{j+\frac{1}{2}} -\nu^{(1)}_{j+\frac{1}{2}}U_{j+\frac{1}{2}}^+ +A(x^-_{j+\frac{1}{2}}) U_{j+\frac{1}{2}}^- -\nu^{(1)}_{j+\frac{1}{2}} U_{j+\frac{1}{2}}^- -\alpha_{1,2}(U_{j+\frac{1}{2}}^+-U_{j+\frac{1}{2}}^-)\right)
\]
where $\alpha_{1,2}=max\{|\lambda^{(1)}(x_{j+\frac{1}{2}})-\nu^{(1)}_{j+\frac{1}{2}}|,|\lambda^{(2)}(x_{j+\frac{1}{2}})-\nu^{(2)}_{j+\frac{1}{2}}|\}$.
Similarly, we can easily update $l^{(2)}_j U_h$ related to $\lambda^{(2)}$ in the following:
\begin{equation}\label{nonmassconserveeqn of v2}
\begin{aligned}
&\frac{d}{dt}\int_{\tilde{I}^{(2)}_j(t)} l^{(2)}_j U_hdx=- \left[l^{(2)}_j \widehat{(A(x)U_h-\nu^{(2)}_{j+\frac{1}{2}}U_h)}|_{j+\frac{1}{2}}\right] + \left[l^{(2)}_{j}\widehat{(A(x)U_h-\nu^{(2)}_{j-\frac{1}{2}}U_h)}|_{j-\frac{1}{2}}\right]\\
  &+ \int_{\tilde{I}^{(2)}_j(t)} l^{(2)}_j F(x,t) dx\doteq L_{2}(U_h(t),t,\tilde{I}^{(2)}_j(t)),
\end{aligned}
\end{equation}
where $l_j^{(2)}$ is a constant vector, locally defined on $\tilde{I}^{(2)}_j(t)$ approximating $l^{(2)}(x)$.
 A simple first order ELDG scheme is composed by two evolution equations \eqref{nonmassconserveeqn of v1} and \eqref{nonmassconserveeqn of v2}. That is, we can update $u_h^1$ by \eqref{nonmassconserveeqn of v1},\eqref{nonmassconserveeqn of v2} and $u_h^1=(r_j^{11}l^{(1)}_j+r_j^{12}l^{(2)}_j)U_h$:
\begin{equation}\label{NMCeqn:integralu1}
 \begin{aligned}
\int_{I_j} u^{1,n+1}_h dx&=\int_{I_j} r_j^{11} l^{(1)}_j U_h^{n+1}dx+\int_{I_j} r_j^{12} l^{(2)}_j U_h^{n+1}dx\\
  &=r_j^{11} \int_{I_j^{*,(1)}} l^{(1)}_j U_h^{n}dx+r_j^{11}\int_{t^n}^{t^{n+1}} L_{1}(U^{(1)}_h(t),t,\tilde{I}^{(1)}_j(t))dt\\
  &+r_j^{12} \int_{I_j^{*,(2)}} l^{(2)}_j U_h^{n}dx+r_j^{12}\int_{t^n}^{t^{n+1}} L_{2}(U^{(2)}_h(t),t,\tilde{I}^{(2)}_j(t))dt,
\end{aligned}
\end{equation}
where $U_h^{n}$ and $U_h^{n+1}$ are defined on the background mesh $I_j$, $U^{(1)}_h(t)$ and $U^{(2)}_h(t)$ are defined on the space-time dynamic meshes $\tilde{I}^{(1)}_j(t)$ and $\tilde{I}^{(2)}_j(t)$ respectively. Similarly, we can update $u_h^2$.

We apply forward-Euler method for time discretization with above ELDG scheme \eqref{NMCeqn:integralu1}:
\begin{equation}
 \begin{aligned}
  &\int_{I_j} u^{1,n+1}dx\\
&=r_j^{11} \int_{I_j^{*,(1)}} l^{(1)}_j U^{n}dx-\Delta t \ r_j^{11}l^{(1)}_j\left[\widehat{(A(x)U^n-\nu^{(1)}_{j+\frac{1}{2}}U^n)}|_{x^{\star, (1)}_{j+\frac12}}-\widehat{(A(x)U-\nu^{(1)}_{j-\frac{1}{2}}U)}|_{x^{\star, (1)}_{j-\frac12}}\right]\\
  &+r_j^{12} \int_{I_j^{*,(2)}} l^{(2)}_j U^{n}dx-\Delta t \ r_j^{12}l^{(2)}_j\left[\widehat{(A(x)U^n-\nu^{(2)}_{j+\frac{1}{2}}U^n)}|_{x^{\star, (2)}_{j+\frac12}}-\widehat{(A(x)U-\nu^{(2)}_{j-\frac{1}{2}}U)}|_{x^{\star, (2)}_{j-\frac12}}\right].\\
\end{aligned}
\end{equation}

\begin{rem}
  The above ELDG scheme is not mass conservative for two reasons:
\begin{itemize}
  \item [(1)] Flux terms can't cancel each other as $r_j^{11}l^{(1)}_j$ and $r_j^{12}l^{(2)}_j$ are discontinuous across cell boundary of $\Omega_j^{(1)}, \Omega_j^{(2)}$.
  \item [(2)] $\sum_j r_j^{11} \int_{I_j^{*,(1)}} l^{(1)}_j U^{n}dx+r_j^{12} \int_{I_j^{*,(2)}} l^{(2)}_j U^{n}dx \neq \sum_j \int_{I_j} U^{n}dx$ because of the inconsistency in characteristic transformations between neighboring cells among $r_j^{11}l^{(1)}_j$ and $r_j^{12}l^{(2)}_j$.
\end{itemize}
\end{rem}

To obtain the mass conservation, a critical point is that eigen-decomposition has to be consistent among two characteristic families and independent of partitions $\Omega_j^{(1)}, \Omega_j^{(2)}$.
%
We choose $a_p(x)$ as a continuous approximation of $a(x)$ which leads that $R_p(x),\Lambda_p(x),R_p^{-1}(x)$ are the approximations of $R(x),\Lambda(x),R^{-1}(x)$ respectively, where
\begin{equation}\label{approxieigenvecterandarepresent}
  R_p(x)=\begin{bmatrix}
r^p_{11}(x)&r^p_{12}(x)\\
r^p_{21}(x)&r^p_{22}(x)
\end{bmatrix}=\begin{bmatrix}
-a_p(x)&a_p(x)\\
1&1
\end{bmatrix}
\end{equation} and
\begin{equation}
  R_p^{-1}(x)=\begin{bmatrix}
l^p_{11}(x)&l^p_{12}(x)\\
l^p_{21}(x)&l^p_{22}(x)
\end{bmatrix}=\begin{bmatrix}
\frac{-1}{2a_p(x)}&\frac12\\
\frac{1}{2a_p(x)}&\frac12
\end{bmatrix} \doteq \begin{bmatrix}
  l_p^{(1)}(x)\\
  l_p^{(2)}(x)
\end{bmatrix}.
\end{equation}
The following equalities hold
\begin{equation}
\begin{aligned}
\label{eq: pickupu1u2}
(r^p_{11}(x)l_p^{(1)}(x)+r^p_{12}(x)l_p^{(2)}(x))U(x)&=u^{1},\\
(r^p_{21}(x)l_p^{(1)}(x)+r^p_{22}(x)l_p^{(2)}(x))U(x)&=u^2,
\end{aligned}
\end{equation}
from $R_p(x)R_p^{-1}(x)=R_p^{-1}(x)R_p(x)=I$. They are critical to design a mass-conservative ELDG scheme.
Take the vector product of $r^p_{11}(x)l_p^{(1)}(x)$ from left with \eqref{eq: 1d_system}, we have a scalar equation
\begin{equation}
  r^p_{11}(x)l_p^{(1)}(x)( U_t+(A(x)U)_x) =r^p_{11}(x)l_p^{(1)}(x) F(x,t).
\end{equation}
Multiply $\psi^{(1)}$ to the equation above,
\begin{equation}\label{eqn: equation1}
  r^p_{11}(x)l_p^{(1)}(x)( U_t+(A(x)U)_x)\psi^{(1)} =r^p_{11}(x)l_p^{(1)}(x) F(x,t)\psi^{(1)}.
\end{equation}
Meanwhile, from $r^p_{11}(x)l_p^{(1)}(x)U$ $\cdot$ (\ref{new ELDG adjointsystem1}), we have
\begin{equation}\label{eqn: equation2}
 r^p_{11}(x)l_p^{(1)}(x)U (\psi^{(1)})_t + r^p_{11}(x)l_p^{(1)}(x)U\alpha^{(1)} (\psi^{(1)})_x =0.
\end{equation}
Next, sum equations (\ref{eqn: equation1}) and (\ref{eqn: equation2}), and integrate over the space-time interval $\Omega^{(1)}_j$, then
\begin{equation}
\begin{aligned}
\label{eq: system_f1}
&\int_{\Omega^{(1)}_j}  \left(r^p_{11}(x)l_p^{(1)}(x)U_t\psi^{(1)}+ r^p_{11}(x)l_p^{(1)}(x)U(\psi^{(1)})_t+ r^p_{11}(x)l_p^{(1)}(x)(A(x)U)_x\psi^{(1)}\right) dxdt\\
&+ \int_{\Omega^{(1)}_j} r^p_{11}(x)l_p^{(1)}(x)U\alpha^{(1)} (\psi^{(1)})_x dxdt=\int_{\Omega^{(1)}_j} r^p_{11}(x)l_p^{(1)}(x) F(x,t)\psi^{(1)} dxdt.
\end{aligned}
\end{equation}
A further manipulation on the L.H.S. of \eqref{eq: system_f1} gives
\begin{equation}
\begin{aligned}
\label{eq:leftsystem_f1}
&\int_{\Omega^{(1)}_j}  \left((r^p_{11}(x)l_p^{(1)}(x)U\psi^{(1)})_t+ r^p_{11}(x)l_p^{(1)}(x)(A(x)U)_x\psi^{(1)}+r^p_{11}(x)l_p^{(1)}(x)U\alpha^{(1)} (\psi^{(1)})_x \right) dxdt\\
=&\int_{\Omega^{(1)}_j}  \left((r^p_{11}(x)l_p^{(1)}(x)U\psi^{(1)})_t+ (r^p_{11}(x)l_p^{(1)}(x)A(x)U\psi^{(1)})_x-(r^p_{11}(x)l_p^{(1)}(x))_xA(x)U\psi^{(1)}\right) dxdt\\
-&\int_{\Omega^{(1)}_j}\left(r^p_{11}(x)l_p^{(1)}(x)A(x)U(\psi^{(1)})_x-r^p_{11}(x)l_p^{(1)}(x)U\alpha^{(1)} (\psi^{(1)})_x\right) dxdt\\
= &\int^{t^{n+1}}_{t^n} \left(\frac{d}{dt}\int_{\tilde{I}^{(1)}_j(t)} (r^p_{11}(x) l_p^{(1)} U\psi^{(1)})dx+[r^p_{11}(x) l_p^{(1)}A(x)U\psi^{(1)}-\nu^{(1)} r^p_{11}(x)l_p^{(1)}U\psi^{(1)}]|_{j-\frac12}^{j+\frac12}\right)dt\\
  &-\int^{t^{n+1}}_{t^n} \int_{\tilde{I}^{(1)}_j(t)}\left((r^p_{11}(x)l_p^{(1)}(x))_xA(x)U\psi^{(1)}+r^p_{11}(x)l_p^{(1)}(x)(A(x)U-\alpha^{(1)}U)\psi^{(1)}_x \right) dxdt.\\
\end{aligned}
\end{equation}
Letting $f^{11}(U)=r^p_{11} l_p^{(1)}(AU-\alpha^{(1)}U)$, the time differential form of \eqref{eq: system_f1} with \eqref{eq:leftsystem_f1} gives
\begin{equation}\label{eqn1 of v1}
\begin{aligned}
&\frac{d}{dt}\int_{\tilde{I}^{(1)}_j(t)} (r^p_{11}(x) l_p^{(1)}(x) U\psi^{(1)})dx+ \left(f^{11}\psi^{(1)}\right) \left|_{\tilde{x}^{(1)}_{j+\frac12}(t) } \right.-\left(f^{11}\psi^{(1)}\right) \left|_{\tilde{x}^{(1)}_{j-\frac12}(t) } \right.
  -\int_{\tilde{I}^{(1)}_j(t)} f^{11}\psi^{(1)}_x dx\\
  &-\int_{\tilde{I}^{(1)}_j(t)}(r^p_{11}(x)l_p^{(1)}(x))_xA(x)U\psi^{(1)} dx= \int_{\tilde{I}^{(1)}_j(t)} r^p_{11}(x)l_p^{(1)}(x) F(x,t)\psi^{(1)} dx.
\end{aligned}
\end{equation}
Similarly, we have an equation related to $\lambda^{(2)}$
\begin{equation}\label{eqn2 of v1:continuous}
\begin{aligned}
&\frac{d}{dt}\int_{\tilde{I}^{(1)}_j(t)} (r^p_{12}(x) l_p^{(2)}(x) U\psi^{(2)})dx+ \left(f^{12}\psi^{(2)}\right) \left|_{\tilde{x}^{(2)}_{j+\frac12}(t) } \right.-\left(f^{12}\psi^{(2)}\right) \left|_{\tilde{x}^{(2)}_{j-\frac12}(t) } \right.
  -\int_{\tilde{I}^{(2)}_j(t)} f^{12}\psi^{(2)}_x dx\\
  &-\int_{\tilde{I}^{(2)}_j(t)}(r^p_{12}(x)l_p^{(2)}(x))_xA(x)U\psi^{(2)} dx= \int_{\tilde{I}^{(2)}_j(t)} r^p_{12}(x)l_p^{(2)}(x) F(x,t)\psi^{(2)} dx,
\end{aligned}
\end{equation}
where $f^{12}(U)=r^p_{12} l_p^{(2)}(AU-\alpha^{(2)}U)$. Then, we can update $u^1$ by \eqref{eq: pickupu1u2},\eqref{eqn1 of v1}, \eqref{eqn2 of v1:continuous}, taking $\Psi^{(1)}(x)=\Psi(x)$ in \eqref{new ELDG adjointsystem1} and $\Psi^{(2)}(x)=\Psi(x)$ in \eqref{new ELDG adjointsystem2}:
\begin{equation}\label{eqn:integralu1continuous}
 \begin{aligned}
  &\int_{I_j} u^{1,n+1}\Psi(x)dx\overset{\eqref{eq: pickupu1u2}}{=}\int_{I_j} r^p_{11} l_p^{(1)}U^{n+1}\Psi(x)dx+\int_{I_j}r^p_{12} l_p^{(2)}U^{n+1}\Psi(x)dx\\
  &=\int_{I_j^{*,(1)}} r^p_{11} l_p^{(1)}U^{n}\psi^{(1)}dx-\int_{t^n}^{t^{n+1}}\left(f^{11}\psi^{(1)}\right) \left|_{\tilde{x}^{(1)}_{j+\frac12}(t) } \right.+\left(f^{11}\psi^{(1)}\right) \left|_{\tilde{x}^{(1)}_{j-\frac12}(t) }\right. dt\\
  &+\int_{t^n}^{t^{n+1}} \int_{\tilde{I}^{(1)}_j(t)} r^p_{11}(x)l_p^{(1)}(x) F(x,t)\psi^{(1)}+(r^p_{11}(x)l_p^{(1)}(x))_xA(x)U\psi^{(1)}+ f^{11}\psi^{(1)}_x dxdt\\
  &+\int_{I_j^{*,(2)}} r^p_{12} l_p^{(2)}U^n\psi^{(2)}dx-\int_{t^n}^{t^{n+1}}\left(f^{12}\psi^{(2)}\right) \left|_{\tilde{x}^{(2)}_{j+\frac12}(t) } \right.+\left(f^{12}\psi^{(2)}\right) \left|_{\tilde{x}^{(2)}_{j-\frac12}(t) }\right. dt\\
  &+\int_{t^n}^{t^{n+1}} \int_{\tilde{I}^{(2)}_j(t)} r^p_{12}(x)l_p^{(2)}(x) F(x,t)\psi^{(2)}+(r^p_{12}(x)l_p^{(2)}(x))_xA(x)U\psi^{(2)}+ f^{12}\psi^{(2)}_x dxdt.
\end{aligned}
\end{equation}
As the ELDG method for scalar, the dynamic interval of $I^{(1)}_j(t)$ and $I^{(2)}_j(t)$ can always be linearly mapped to a reference cell $\xi\in I_j$, the ELDG discretization of eq.~\eqref{eqn:integralu1continuous} is to find $u_h^1(x,t)\in P^k(I_j(t))$, so that
\begin{equation}\label{eqn:integralu1}
 \begin{aligned}
  \int_{I_j} u^{1}_h(x,t^{n+1})\Psi(x)dx&=\int_{I_j^{*,(1)}} r^p_{11} l_p^{(1)}U^{n}_h\psi^{(1)}dx+\int_{t^n}^{t^{n+1}} L_{11}(U_h(t),t,\tilde{I}^{(1)}_j(t))dt\\
  &+\int_{I_j^{*,(2)}} r^p_{12} l_p^{(2)}U^n_h\psi^{(2)}dx+\int_{t^n}^{t^{n+1}} L_{12}(U_h(t),t,\tilde{I}^{(2)}_j(t))dt,
\end{aligned}
\end{equation}
for $\psi^{(1)}(x, t)$ satisfying the adjoint problem \eqref{new ELDG adjointsystem1} with $\forall \Psi(x) = \psi(x, t^{n+1})\in P^k(I_j)$. Here
\begin{equation}
 \begin{aligned}
  L_{11}(U_h(t),t,\tilde{I}^{(1)}_j(t))&=-\widehat{f^{11}_{j+\frac{1}{2}}} \psi^{(1),-}_{j+\frac{1}{2}}+\widehat{f^{11}_{j-\frac{1}{2}}} \psi^{(1),+}_{j-\frac{1}{2}}+\int_{\tilde{I}^{(1)}_j(t)}f^{11}\psi^{(1)}_{x}(x,t)dx\\
  &+\int_{\tilde{I}^{(1)}_j(t)}(r^p_{11} l_p^{(1)})_{x} AU_h\psi^{(1)}(x,t)+r^p_{11} l_p^{(1)} F\psi^{(1)}(x,t) dx,\\
 L_{12}(U_h(t),t,\tilde{I}^{(2)}_j(t))&=-\widehat{f^{12}_{j+\frac{1}{2}}} \psi^{(2),-}_{j+\frac{1}{2}}+\widehat{f^{12}_{j-\frac{1}{2}}} \psi^{(2),+}_{j-\frac{1}{2}}+\int_{\tilde{I}^{(2)}_j(t)}f^{12}\psi^{(2)}_{x}(x,t)dx\\
  &+\int_{\tilde{I}^{(2)}_j(t)}(r^p_{12} l_p^{(2)})_{x} AU_h\psi^{(2)}(x,t)+r^p_{12} l_p^{(2)} F\psi^{(2)}(x,t) dx,
\end{aligned}
\end{equation}
where
 $\psi^{(1),\pm}_{j\pm\frac{1}{2}}=\psi^{(1)}(x^{(1),\pm}_{j\pm\frac{1}{2}}(t),t)$ and $\widehat{f^{11}}$ at a cell boundary can be taken as a monotone flux, e.g. the Lax-Friedrichs flux
\[
\widehat{f^{11}_{j+\frac{1}{2}}}(U)=r^p_{11}(x_{j+\frac{1}{2}})l_p^{(1)}(x_{j+\frac{1}{2}})
\widehat{(A\cdot U-\nu^{(1)}_{j+\frac{1}{2}}U)}_{{j+\frac{1}{2}}}.
\]
$\psi^{(2),\pm}_{j\pm\frac{1}{2}}=\psi^{(2)}(x^{(2),\pm}_{j\pm\frac{1}{2}}(t),t)$ and $\widehat{f^{12}}$ can be similarly defined at a cell boundary.
We can similarly obtain the ELDG scheme for $u^{2}_h$.

\noindent {\bf (4) Fully discrete ELDG scheme with method-of-lines RK schemes.}

 To update \eqref{eqn:integralu1} from $U^{n}_h$ to $U^{n+1}_h$, we first apply the forward Euler time discretization to get 1st order accuracy, then we generalize the scheme to general RK methods.
 There are two main steps involved here.
 In order to describe the implementation procedure of the fully discrete ELDG scheme, we define the $L^2$ projection.
 \begin{defn}($L^2$ projection)
   Let $u(x)\in L^2(\Omega)$, $M=\{ I_j \}_{j=1}^N$ and $\tilde{M}=\{ \tilde{I}_j \}_{j=1}^N$ be two different meshes of $\Omega$. We have function spaces $V_h^k=\{u(x): u|_{I_j}\in P^k(I_j), \forall j \}$ and $\tilde{V}_h^k=\{\tilde{u}(x): \tilde{u}|_{\tilde{I}_j}\in P^k(\tilde{I}_j), \forall j \}$ corresponding to meshes $M$ and $\tilde{M}$. The $L^2$ projection of $u_M\in V_h^k$ onto space $\tilde{V}_h^k$ can be defined as, find $\tilde{u}_{\tilde{M}}(x) \in\tilde{V}_h^k$, s.t.
   \beq\label{eqn:L2projectiondefination}
   \int_{\tilde{I}_j} \tilde{u}_{\tilde{M}}(x) \varphi(x) dx=\int_{\tilde{I}_j} u_M(x) \varphi(x) dx, \ \forall \varphi(x)\in \tilde{V}_h^k.
   \eeq
 We denote $\tilde{u}_{\tilde{M}}(x)=Proj[ u_M(x); M, \tilde{M} ]$. The evaluation of the right hand side of ~\eqref{eqn:L2projectiondefination} can be done in a subinterval-by-subinterval fashion. The implementation details can be found in \cite{Guo2013discontinuous}.
 \end{defn}
 Then, we propose a fully discrete EL RK DG scheme with procedure as described:
\begin{enumerate}
\item {\bf Obtain the initial condition} on upstream meshes ${\tilde{I}^{(1)}_j(t^n)}$ and ${\tilde{I}^{(2)}_j(t^n)}$ of \eqref{eqn:integralu1} by $U^{(1)}_h(t^n)=Proj[ U^{n}_h; I_j, {\tilde{I}^{(1)}_j(t^n)}]$ and $U^{(2)}_h(t^n)=Proj[ U^{n}_h; I_j, {\tilde{I}^{(2)}_j(t^n)}]$, which are the $L^2$ projections of solutions from the background mesh to the upstream mesh.
\item {\bf Update \eqref{eqn:integralu1} from $U^{n}_h$ to $U^{n+1}_h$, component-by-component.}
\begin{enumerate}
  \item[(a)] Get the mesh information of the dynamic element $\tilde{I}^{(1)}_j(t^{(l)}), l=0,...,s$ on RK stages by $\tilde{x}^{(1)}_{j\pm\frac12}(t) = x_{j\pm\frac12} + (t-t^{n+1}) \nu^{(1)}_{j\pm\frac12}$. Here $s=1$ for forward-Euler method and $s=2$ for SSPRK2, see the blue domain in Figure~\ref{iso2} and for general RK2 with intermediate stage in Figure~\ref{iso16}.
  \item[(b1)] For forward-Euler method, compute
  \begin{equation}\label{eq: forward-euler for system}
 \begin{aligned}
  &\int_{I_j} u^{1,n+1}_h\Psi(x)dx=\int_{I_j^{*,(1)}} r^p_{11} l_p^{(1)}U^{n}_h\psi^{(1),n}dx+\Delta t L_{11}(U_h^n,t^n,\tilde{I}^{(1)}_j(t^n))\\
  +&\int_{I_j^{*,(2)}} r^p_{12} l_p^{(2)}U^{n}_h\psi^{(2),n}dx+\Delta t L_{12}(U_h^n,t^n,\tilde{I}^{(2)}_j(t^n))\\
  =&{\color{blue}\int_{I_j^{*,(1)}} r^p_{11} l_p^{(1)}U^{n}_h\psi^{(1),n}dx}+\Delta t \left(
  -\widehat{f^{11}_{j+\frac{1}{2}}}(U_h^{(1)}(t^n)) \psi^{(1),n,-}_{j+\frac{1}{2}}+\widehat{f^{11}_{j-\frac{1}{2}}}(U_h^{(1)}(t^n)) \psi^{(1),n,+}_{j-\frac{1}{2}} \right.\\
  &\left. +\int_{I_j^{*,(1)}} f^{11}(U_h^{(1)}(t^n))\psi^{(1)}_{x}(x,t^n)dx+{\color{blue}\int_{I_j^{*,(1)}} (r^p_{11} l_p^{(1)})_{x} AU^n_h\psi^{(1)}(x,t^n) dx}\right.\\
  &\left. + \int_{I_j^{*,(1)}} r^p_{11} l_p^{(1)} F(x,t^n)\psi^{(1)}(x,t^n) dx\right)\\
  &+{\color{blue}\int_{I_j^{*,(2)}} r^p_{12} l_p^{(2)}U^{n}_h\psi^{(2),n}dx}+\Delta t \left(
  -\widehat{f^{12}_{j+\frac{1}{2}}}(U_h^{(2)}(t^n)) \psi^{(2),n,-}_{j+\frac{1}{2}}+\widehat{f^{12}_{j-\frac{1}{2}}}(U_h^{(2)}(t^n)) \psi^{(2),n,+}_{j-\frac{1}{2}} \right.\\
  &\left. +\int_{I_j^{*,(2)}} f^{12}(U_h^{(2)}(t^n))\psi^{(2)}_{x}(x,t^n)dx+{\color{blue}\int_{I_j^{*,(2)}} (r^p_{12} l_p^{(2)})_{x} AU^n_h\psi^{(2)}(x,t^n) dx}\right.\\
  &\left. + \int_{I_j^{*,(2)}} r^p_{12} l_p^{(2)} F(x,t^n)\psi^{(2)}(x,t^n) dx\right)\\
\end{aligned}
\end{equation}
where $\Delta t^n=t^{n+1}-t^n$, $U_h^{(1)}(t^n)=Proj[U^{n}_h,\{ I_j \}_{j=1}^N, \{I_j^{*,(1)}\}_{j=1}^N]$ and \\
$U_h^{(2)}(t^n)=Proj[U^{n}_h,\{ I_j \}_{j=1}^N, \{I_j^{*,(2)}\}_{j=1}^N]$.
We compute the four integration terms of ~\eqref{eq: forward-euler for system}$\int_{I_j^{*,(1)}} r^p_{11} l_p^{(1)}U^{n}_h\psi^{(1),n}dx$, $\int_{I_j^{*,(2)}} r^p_{12} l_p^{(2)}U^{n}_h\psi^{(2),n}dx$, $\int_{\tilde{I}^{(1)}_j(t^n)}(r^p_{11} l_p^{(1)})_{x} AU^n_h\psi^{(1)}(x,t^n) dx$ and $\int_{\tilde{I}^{(2)}_j(t^n)}(r^p_{12} l_p^{(2)})_{x} AU^n_h\psi^{(2)}(x,t^n) dx$ highlighted in blue with $U_h^n$ on background meshes in order to makes the scheme satisfy mass conservation.
Thus they have to be evaluated subinterval-by-subinterval since DG solution is discontinuous.
In summary, we can get $u_h^{1,n+1}$ on $I_j$ by \eqref{eq: forward-euler for system}. Similarly, we can get $u_h^{2,n+1}$ on $I_j$.

\begin{minipage}{\textwidth}
\begin{minipage}[t]{0.48\textwidth}
\makeatletter\def\@captype{table}
	\begin{tabular}{c|c c}
		0 & & \\
		1 & 1 & \\
		\hline
		& $\frac{1}{2}$ & $\frac{1}{2}$
	\end{tabular}
\caption{SSPRK2 Butcher Tableau}
\label{Table:SSPRK2ButcherTable}
\end{minipage}
\begin{minipage}[t]{0.48\textwidth}
	\makeatletter\def\@captype{table}
	\begin{tabular}{c|c c}
		0 & & \\
		$\frac{1}{2}$ & $\frac{1}{2}$ & \\
		\hline
		& 0 & 1
	\end{tabular}
	\caption{RK2 Butcher Tableau}
	\label{Table:RK2ButcherTable}
\end{minipage}
\end{minipage}

\item[(b2)] For SSPRK2 method with Butcher tableau: \ref{Table:SSPRK2ButcherTable}, we get $u^{1}_h(t^{(1)})$ from \eqref{eq: forward-euler for system}, then compute
\begin{equation}
\begin{aligned}
  &\int_{I_j} u^{1,n+1}_h\Psi(x)dx=\int_{I_j^{*,(1)}} r^p_{11} l_p^{(1)}U^{n}_h\psi^{(1),n}dx+0.5\Delta t L_{11}(U_h(t^{(1)}),t^{(1)},I_j)\\
  +&\int_{I_j^{*,(2)}} r^p_{12} l_p^{(2)}U^{n}_h\psi^{(2),n}dx+0.5\Delta t L_{12}(U_h(t^{(1)}),t^{(1)},I_j),\\
\end{aligned}\label{eq: sspRK2_system}
\end{equation}
where 
$t^{(1)}=t^{n+1}$, $u^{1}_h(t^{(1)})$ and $u^{2}_h(t^{(1)})$ are defined on background mesh $I_j$.
\item[(b3)] For general RK methods with intermediate stages, we will update intermediate RK solutions on background mesh as in \cite{ding2020semi}. For example, for a 2nd order RK method with Butcher tableau \ref{Table:RK2ButcherTable}. It has an intermediate stage at $t^{(1)}=t^{n+\frac{\Delta t}{2}}$, we propose the following steps as \cite{ding2020semi}, see Figure ~\ref{iso16}.
\begin{enumerate}
  \item We denote the dynamic domain tracking $I_j$ from $t^{(1)}$ to $t_n$ with speed $\nu^{(1)}_{j\pm \frac12}$ at mesh point $x_{j\pm \frac12}$ as $\tilde{I}^{(1)}_{j,(1)}(t)$, see the green domain in Figure~\ref{iso16}. $\tilde{I}^{(2)}_{j,(1)}(t)$ related to second characteristic is defined similarly. 
      Then we can update $U^{(1)}_h$ on $\tilde{I}^{(1)}_{j,(1)}(t)$ and $\tilde{I}^{(2)}_{j,(1)}(t)$ from $t^n$ to $t^{(1)}$ as in a forward-Euler method.
  \item We update $U_h^{n+1}$ on dynamic domain $\tilde{I}^{(1)}_j(t)$ and $\tilde{I}^{(2)}_j(t)$ from $U_h^{n}$ with projection onto $I_j^{*,(1)}$ and $U_h(t^{(1)})$ with projection onto $\tilde{I}^{(1)}_{j,(1)}(t^{(1)})$.
\end{enumerate}
  \end{enumerate}
\end{enumerate}

\begin{figure}[h!]
\centering
\begin{tikzpicture}
[x=1cm,y=1cm]
  \begin{scope}[thick]

  \draw[fill=blue!5] (0.,2) -- (-1.5,0) -- (1.6,0) -- (2.4,2)
      -- cycle;
  \draw[fill=green!5] (0.,1) -- (-0.75,0) -- (2.0,0) -- (2.4,1)
      -- cycle;
   \node[green!70, rotate=0] (a) at ( 1. ,0.5) { ${ \tilde{I}^{(1)}_{j,(1)} }(t)$ };

   \draw (-3,3) node[fill=white] {};
    \draw (-3,-1) node[fill=white] {};
    \draw[black]                   (-3,0) node[left] {} -- (3,0)
                                        node[right]{$t^{n}$};
    \draw[black] (-3,2) node[left] {$$} -- (3,2)
                                        node[right]{$t^{n+1}$};

     \draw[snake=ticks,segment length=2.4cm] (-2.4,2) -- (0,2) node[above left] {$x_{j-\frac12}$};
     \draw[snake=ticks,segment length=2.4cm] (0,2) -- (2.4,2) node[above right] {$x_{j+\frac12}$};

          \draw[snake=ticks,segment length=2.4cm] (-2.4,0) -- (0,0);
     \draw[snake=ticks,segment length=2.4cm] (0,0) -- (2.4,0);

            \draw[blue,thick] (0.,2) node[left] {$$} -- (-1.5,0)
                                        node[black,below]{$x^{\star,(1)}_{j-\frac12}$ };

                \draw[blue,thick] (2.4,2) node[left] {$$} -- (1.6,0)
                                        node[black,below right]{ $x^{\star,(1)}_{j+\frac12}$};
\draw [decorate,color=blue,decoration={brace,mirror,amplitude=9pt},xshift=0pt,yshift=0pt]
(-1.5,0) -- (1.6,0) node [blue,midway,xshift=0cm,yshift=-20pt]
{\footnotesize $\tilde{I}^{(1)}_j(t^{n}) = I^{\star,(1)}_{j}$};

\draw [decorate,color=blue,decoration={brace,amplitude=10pt},xshift=0pt,yshift=0pt]
(0,2) -- (2.4,2) node [blue,midway,xshift=0cm,yshift=20pt]
{\footnotesize $\tilde{I}^{(1)}_j(t^{n+1}) = I_{j}$};

\draw [decorate,color=blue,decoration={brace,amplitude=10pt},xshift=0pt,yshift=0pt]
(-0.75,1) -- (2.,1) node [blue,midway,xshift=0cm,yshift=20pt]
{\footnotesize $\tilde{I}^{(1)}_j(t^{(1)})$};

%

  \node[blue, rotate=54] (a) at (-1.0,1.2) { $\nu^{(1)}_{j-\frac12}$ };
  \node[blue, rotate=70] (a) at (2.25,1.) { $\nu^{(1)}_{j+\frac12}$ };

  \draw[dashed,gray] (0.,2) -- (0,0);
  \draw[dashed,gray] (2.4,2) -- (2.4,0);

\draw[blue,thick] (-0.75,1) -- (2.0,1);
\draw[red,thick] (0.,1) -- (2.4,1);
\draw[green,thick] (0.,1) -- (-0.75,0);
\draw[green,thick] (2.4,1) -- (2.0,0);

  \end{scope}
\end{tikzpicture}
\caption{Update RK intermediate solution at background mesh (red line) from the first characteristic family of a hyperbolic system.}
\label{iso16}
\end{figure}
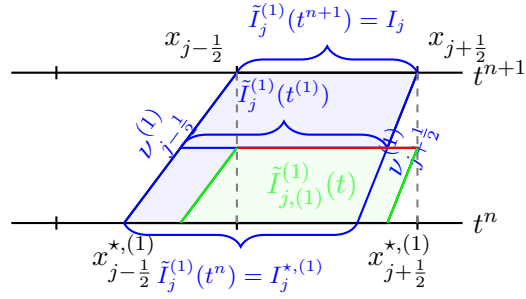
\begin{thm} (Mass conservation)
The proposed fully discrete ELDG scheme is locally mass conservative.
In particular, given a DG solution $u_h(x,t^n)\in V_h^k$ with a periodic boundary condition, we have
\begin{equation*}
\sum_{ i=1 }^N \int_{ I_j } U_h( x,t^{n+1} ) dx
=
\sum_{ i=1 }^N \int_{ I_j } U_h( x,t^n ) dx.
\end{equation*}
\end{thm}
\noindent
{\em Proof.} We firstly consider the forward Euler time discretization. Taking $\Psi=1$ and $F=0$ in the scheme \eqref{eq: forward-euler for system}, we have
\begin{equation}
 \begin{aligned}
  &\int_{\Omega} u^{1,n+1}_hdx=\sum_j\int_{I_j} u^{1,n+1}_hdx\\
  &=\sum_j\left[\int_{I_j^{*,(1)}} r^p_{11} l_p^{(1)}U^{n}_hdx+\Delta t\left( -\widehat{f^{11}_{j+\frac{1}{2}}}(U_h^{(1)}(t^n))+\widehat{f^{11}_{j-\frac{1}{2}}}(U_h^{(1)}(t^n))+\int_{I_j^{*,(1)}} (r^p_{11} l_p^{(1)})_{x} AU_h^n dx \right)\right]\\
  &+\sum_j\left[\int_{I_j^{*,(2)}} r^p_{12} l_p^{(2)}U^{n}_hdx+\Delta t\left( -\widehat{f^{12}_{j+\frac{1}{2}}}(U_h^{(2)}(t^n))+\widehat{f^{12}_{j-\frac{1}{2}}}(U_h^{(2)}(t^n))+\int_{I_j^{*,(2)}} (r^p_{12} l_p^{(2)})_{x} AU^n_h dx \right)\right]\\
  &=\sum_j\left(\int_{I_j^{*,(1)}} r^p_{11} l_p^{(1)}U^{n}_hdx+\int_{I_j^{*,(2)}} r^p_{12} l_p^{(2)}U^{n}_hdx\right)\\
  &+\Delta t\sum_j\left(\int_{I_j^{*,(1)}} (r^p_{11} l_p^{(1)})_{x} AU_h^n dx+\int_{I_j^{*,(2)}} (r^p_{12} l_p^{(2)})_{x} AU^n_h dx\right)\\
  &=\int_{\Omega} (r^p_{11}l_p^{(1)}+r^p_{12} l_p^{(2)})U^{n}_h dx+\Delta t\int_{\Omega} (r^p_{11} l_p^{(1)}+r^p_{12} l_p^{(2)})_{x} AU_h^n dx\\
  &=\int_{\Omega} u^{1,n}_hdx,
\end{aligned}
\end{equation}
which follows from the cancellation of unique fluxes at cell boundaries, $r^p_{11}l_p^{(1)}+r^p_{12} l_p^{(2)}=[1,0]$ and \eqref{eq: pickupu1u2} with integration in a subinterval-by-subinterval fashion.
The mass conservation for the fully discrete ELDG scheme can be proved in a similar fashion.

\begin{rem}
   To maintain the mass-conservative property, the choice of eigenvectors $R_p(x)$ is not necessarily exact for ELDG scheme, as long as $R_p(x)$ and $R_p^{-1}(x)$ are a consistent pair throughout the domain.
\end{rem}

\subsection{2D linear wave equation}
The solution for high-dimensional hyperbolic systems is given by means of a characteristic cone, rather than individual characteristic lines.
Numerically, characteristic Galerkin \cite{ostkamp1997multidimensional} or evolution Galerkin \cite{butler1960numerical,lax1959systems,lukavcova2000evolution} methods have been proposed and developed to solve high dimensional hyperbolic system. This method is constructed by taking into account information propagated in all bicharacteristic directions and involving integrals around the characteristic cone. However, the awkward integrals over the mantle, involving intermediate times, limit both the accuracy and the stability of the resulting schemes. Thus the finite volume evolution Galerkin (FVEG) schemes are introduced, which is in a predictor-corrector plus finite volume framework to get higher accuracy. FVEG method has been developed in \cite{lukavcova2007finite} and widely applied in \cite{kroger2005evolution,lukavcova2007finite,lukavcova2007well}, and the stability and accuracy have been investigated in \cite{lukavcova2002finite,lukavcova2007finite,lukavcova2006stability,lukacova2009entropy}. Even though the FVEG method can achieve high-order accuracy and stability with extra large step, the algorithm implementation is still very complex for high-dimensional problems. In this paper, we use the splitting method for higher dimensional problem.

Consider a linear 2D wave equation
\begin{equation}\label{eqn: 2d linear wave equation}
  u_{tt}= (a(x,y,t)u_x)_x + (b(x,y,t)u_y)_y, (x,y) \in \Omega.
\end{equation}
It can be written as a first order linear hyperbolic system
\begin{equation}\label{eqn: 2d linear equation}
  U_t+(A(x,y,t)U)_x + (B(x,y,t)U)_y=0, (x,y) \in \Omega,
\end{equation}
where \[
U=\begin{bmatrix}
  u_t\\
  u_x\\
  u_y
\end{bmatrix},
\qquad
A=\begin{bmatrix}
0&-a(x,y,t)&0\\
-1&0&0\\
0&0&0
\end{bmatrix},
\qquad
B=\begin{bmatrix}
0&0&-b(x,y,t)\\
0&0&0\\
-1&0&0
\end{bmatrix}.
\]
The domain $\Omega$ is partitioned into rectangular meshes with each computational cell $\Omega_{ij} = [x_{i- \frac{1}{2}}, x_{i+\frac{1}{2}}] \times [y_{j+ \frac{1}{2}}, y_{j+\frac{1}{2}}]$, where we use the piecewise $Q^k$ tensor-product polynomial spaces. Then we extend ELDG algorithm to 2D problems via dimensional splitting \cite{qiu2011positivity}.
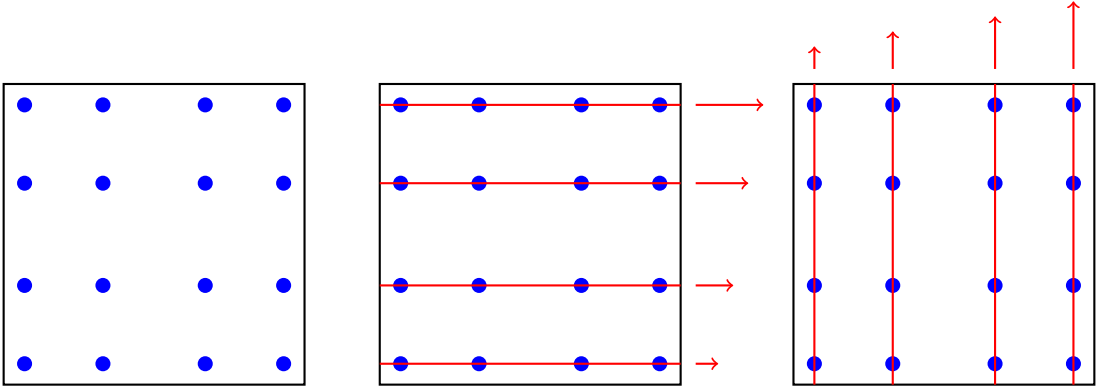
\begin{figure}[h!]
\centering
\begin{tikzpicture}[x=1cm,y=1cm]
  \begin{scope}[thick]

  \draw[fill=white!5] (-2.,-2.) -- (-2.,2) -- (2.,2.) -- (2.,-2.)
      -- cycle;
  \fill[blue] (-1.722,-1.722) circle (0.1);
  \fill[blue] (-1.722,-0.68) circle (0.1);
  \fill[blue] (-1.722,1.722) circle (0.1);
  \fill[blue] (-1.722,0.68) circle (0.1);
  \fill[blue] (-0.68,-1.722) circle (0.1);
  \fill[blue] (-0.68,-0.68) circle (0.1);
  \fill[blue] (-0.68,0.68) circle (0.1);
  \fill[blue] (-0.68,1.722) circle (0.1);
  \fill[blue] (0.68,-1.722) circle (0.1);
  \fill[blue] (0.68,-0.68) circle (0.1);
  \fill[blue] (0.68,0.68) circle (0.1);
  \fill[blue] (0.68,1.722) circle (0.1);
  \fill[blue] (1.722,-1.722) circle (0.1);
  \fill[blue] (1.722,-0.68) circle (0.1);
  \fill[blue] (1.722,1.722) circle (0.1);
  \fill[blue] (1.722,0.68) circle (0.1);
  \end{scope}

\begin{scope}[thick]
 \draw[fill=white!5] (3.,-2.) -- (3.,2) -- (7.,2.) -- (7.,-2.)
      -- cycle;
  \fill[blue] (3.278,-1.722) circle (0.1);
  \fill[blue] (3.278,-0.68) circle (0.1);
  \fill[blue] (3.278,1.722) circle (0.1);
  \fill[blue] (3.278,0.68) circle (0.1);
  \fill[blue] (4.32,-1.722) circle (0.1);
  \fill[blue] (4.32,-0.68) circle (0.1);
  \fill[blue] (4.32,0.68) circle (0.1);
  \fill[blue] (4.32,1.722) circle (0.1);
  \fill[blue] (5.68,-1.722) circle (0.1);
  \fill[blue] (5.68,-0.68) circle (0.1);
  \fill[blue] (5.68,0.68) circle (0.1);
  \fill[blue] (5.68,1.722) circle (0.1);
  \fill[blue] (6.722,-1.722) circle (0.1);
  \fill[blue] (6.722,-0.68) circle (0.1);
  \fill[blue] (6.722,1.722) circle (0.1);
  \fill[blue] (6.722,0.68) circle (0.1);
  \draw[red] (3.,-1.722) -- (7.,-1.722);
  \draw[red] (3.,-0.68) -- (7.,-0.68);
  \draw[red] (3.,0.68) -- (7.,0.68);
  \draw[red] (3.,1.722) -- (7.,1.722);
  \draw[red][->] (7.2,-1.722) -- (7.5,-1.722);
  \draw[red][->] (7.2,-0.68) -- (7.7,-0.68);
  \draw[red][->] (7.2,0.68) -- (7.9,0.68);
  \draw[red][->] (7.2,1.722) -- (8.1,1.722);
\end{scope}
\begin{scope}[thick]
 \draw[fill=white!5] (8.5,-2.) -- (8.5,2) -- (12.5,2.) -- (12.5,-2.)
      -- cycle;
  \fill[blue] (8.778,-1.722) circle (0.1);
  \fill[blue] (8.778,-0.68) circle (0.1);
  \fill[blue] (8.778,1.722) circle (0.1);
  \fill[blue] (8.778,0.68) circle (0.1);
  \fill[blue] (9.82,-1.722) circle (0.1);
  \fill[blue] (9.82,-0.68) circle (0.1);
  \fill[blue] (9.82,0.68) circle (0.1);
  \fill[blue] (9.82,1.722) circle (0.1);
  \fill[blue] (11.18,-1.722) circle (0.1);
  \fill[blue] (11.18,-0.68) circle (0.1);
  \fill[blue] (11.18,0.68) circle (0.1);
  \fill[blue] (11.18,1.722) circle (0.1);
  \fill[blue] (12.222,-1.722) circle (0.1);
  \fill[blue] (12.222,-0.68) circle (0.1);
  \fill[blue] (12.222,1.722) circle (0.1);
  \fill[blue] (12.222,0.68) circle (0.1);
  \draw[red] (8.778,-2.) -- (8.778,2.);
  \draw[red] (9.82,-2.) -- (9.82,2.);
  \draw[red] (11.18,-2.) -- (11.18,2.);
  \draw[red] (12.222,-2.) -- (12.222,2.);
  \draw[red][->] (8.778,2.2) -- (8.778,2.5);
  \draw[red][->] (9.82,2.2) -- (9.82,2.7);
  \draw[red][->] (11.18,2.2) -- (11.18,2.9);
  \draw[red][->] (12.222,2.2) -- (12.222,3.1);

\end{scope}
\end{tikzpicture}
\caption{Illustration of the 2D ELDG scheme via Strang splitting. k = 3. }
\label{splittingk3}
\end{figure}
\begin{enumerate}
 \item We first locate $(k + 1)^2$ tensor-product Gaussian nodes on cell $\Omega_{ij}: (x_{i,p}, y_{j,q}),\  p,q = 0,... ,k$. For example, see Figure \ref{splittingk3} (left) for the case of $k=3$.
 \item Then, the equation \eqref{eqn: 2d linear equation} is split into two 1D hyperbolic problems based on the quadrature nodes in both $x-$ and $y-$ directions:
\begin{align}\label{eqn: split the 2d linear equation 1}
  U_t+(A(x,y,t)U)_x &= 0,\\
 \label{eqn: split the 2d linear equation 2}
  U_t+(B(x,y,t)U)_y &= 0.
\end{align}
Based on a 1D ELDG formulation, the split equations \eqref{eqn: split the 2d linear equation 1} and \eqref{eqn: split the 2d linear equation 2} are evolved via Strang splitting
over a time step $\Delta t$ as follows.
\begin{itemize}
  \item Evolve 1D equation \eqref{eqn: split the 2d linear equation 1} at different $y_{j,q}'s$ for a half time-step $\Delta t/2$, see Figure \ref{splittingk3} (middle). For each $y_{j,q}$, the $(k + 1)$ point values are mapped to a $P^k$ polynomial per cell, then the 1D system \eqref{eqn: split the 2d linear equation 1} is evolved by the proposed ELDG scheme. Finally, we can map the evolved $P^k$ polynomial back to the $(k + 1)$ point values to update the solution.
  \item Evolve 1D system \eqref{eqn: split the 2d linear equation 2} at different $x_{i,p}'s$ for a full time-step $\Delta t$ as above, see Figure \ref{splittingk3} (right).
  \item Evolve 1D system \eqref{eqn: split the 2d linear equation 1} at different $y_{j,q}'s$ for another half time-step $\Delta t/2$.
\end{itemize}
\end{enumerate}
The splitting 2D ELDG formulation maintains
high order accuracy in space, extra large time stepping size with stability and mass conservation; and has a second order splitting error.

\section{Numerical results}
\label{section:numerical}

In this section, we show numerical results of the proposed scheme for wave equations. We set the time stepping size as $\Delta t= \frac{CFL}{a} \Delta x$ for 1D and $\Delta t=\frac{CFL}{\frac{a}{\Delta x}+\frac{b}{\Delta y}}$, where a and b are maximum eigenvalues of coefficient matrixes in $x$- and $y$-directions respectively.
We use RK4 \cite{Jiang_Shu} for time discretization.
We study the following aspects: the spatial order of convergence by using small enough time stepping sizes, the spatial super-convergence of the post-processed solutions \cite{bramble1977higher,cockburn2003enhanced} produced by convolving the ELDG solution with a suitable kernel consisting of B-splines for the purpose of improving spatial convergence rate, the temporal order of convergence and numerical stability under a large time stepping size by varying CFL for a fixed spatial mesh.
%

\subsection{1D wave equations}
\begin{exa}
(1D wave equation with constant coefficient.)
We consider the 1D wave equation eq.~\eqref{2d_hcl} with constant coefficient $a(x)=1$ and the source term $f(x,t)=0$. The initial data is $u(x,0)=\sin(x), x\in [0,2\pi]$ with periodic boundary condition. The exact solution is $u(x,t)=\sin(x+t)$.
For the constant coefficient problem, if using exact characteristic velocity fields for space-time partition and exact eigenvectors, the proposed ELDG method is the same as SL DG, then it is unconditionally stable.
Here we perturb the characteristic velocity $\nu_{j+\frac12}^{(1)}$ in \eqref{eq: ELDG alpha for system_def1} at cell boundaries and/or $a_p(x)$ in \eqref{approxieigenvecterandarepresent} related to approximating eigenvectors to get ELDG, ELDG1, ELDG2 and ELDG3 schemes respectively. Similarly we implement the non mass-conservative ELDG methods denoted as NMC ELDG, NMC ELDG1, NMC ELDG2 and NMC ELDG3.
Related parameters of these ELDG methods are given in Table \ref{The related parameters settings for system}.
\begin{table}[!ht]\footnotesize
\caption{The numerical parameters of ELDG, ELDG1, ELDG2, ELDG3 methods and NMC ELDG, NMC ELDG1, NMC ELDG2, NMC ELDG3 methods for $u_{tt}=u_{xx}$.
  }
\centering
\begin{tabular}{| c | c | c| c| c| c| }

 \hline
  setting  &  \multicolumn{1}{c|}{ ELDG} & \multicolumn{1}{c|}{ ELDG1}& \multicolumn{1}{c|}{ ELDG2}& \multicolumn{1}{c|}{ ELDG3}
  \\ \hline
    $\nu^{(1)}_{j+\frac12}$ in \eqref{eq: ELDG alpha for system_def1}, $\nu^{(2)}_{j+\frac12}=-\nu^{(1)}_{j+\frac12}$ &  $1$   &  $1+\sin(x_j)\Delta x$    &  $1$   &  $1+\sin(x_j)\Delta x$\\
  \hline
    $a_p(x)$ in \eqref{approxieigenvecterandarepresent}  &    $1$  &    $1$     &     $1+\sin(x)\Delta x$    &    $1+\sin(x)\Delta x$  \\
 \hline
  setting  &  \multicolumn{1}{c|}{ NMC ELDG} & \multicolumn{1}{c|}{NMC ELDG1}& \multicolumn{1}{c|}{NMC ELDG2}& \multicolumn{1}{c|}{NMC ELDG3}\\
  \hline
    $\nu^{(1)}_{j+\frac12}$ in \eqref{eq: ELDG alpha for system_def1}, $\nu^{(2)}_{j+\frac12}=-\nu^{(1)}_{j+\frac12}$ &  $1$   &  $1+\sin(x_j)\Delta x$   &  $1$   &  $1+\sin(x_j)\Delta x$\\
  \hline
   $a_j$ in \eqref{nmcapproxieigenvectora}  &    $1$  &    $1$     &     $1+\sin(x_j)\Delta x$    &    $1+\sin(x_j)\Delta x$  \\
\hline
  \end{tabular}
\label{The related parameters settings for system}
\end{table}

Table \ref{linear1d_spatial system for 'background scheme'} and \ref{linear1d_spatial system for 'background scheme'post} report spatial accuracies of the ELDG, ELDG1, ELDG2 and ELDG3 methods for this example under the same time stepping size without and with post-processing technique \cite{bramble1977higher,cockburn2003enhanced}. We can observe the optimal convergence rate $k+1$ and $2k+1$.
We vary time stepping sizes, with fixed well-resolved spatial meshes, and plot error vs. $CFL$ in Figure \ref{linear_stability1 for system} and \ref{linear_stability1 for systempost} for ELDG, ELDG1, ELDG2 and ELDG3 schemes without and with post-processed technique respectively, after a long time $T=100$. 
The plots from post-processed ELDG schemes better show the fourth order temporal convergence. ELDG2 and ELDG3 perform comparably; they have a more restricted time step constraint than ELDG1.
It indicates that, stability is affected by approximations of characteristic via the space-time partition and approximation of eigenvectors.
We also note that, in both Figure \ref{linear_stability1 for system} and \ref{linear_stability1 for systempost}, the CFL allowed with stability is much larger than that of the RK DG method which is $\frac{1}{2k+1}$. Further, we verify the mass conservative property of the ELDG schemes are around machine precision and the non-mass conservative property of the NMC ELDG schemes is presented in Figure \ref{linearmassconservesmoothinitial}.
\begin{table}[!ht]\footnotesize
\caption{1D wave equation with constant coefficient. $ u_{tt}  =  u_{xx} $ with initial condition $u(x,0) = \sin(x)$. $T=1$.
We use $CFL=0.3$ and $CFL=0.18$ with RK4 time discretization for all $P^1$ and $P^2$ respectively. The error for only $u^1=u_t$ was shown in this table.
  }
\centering
\begin{tabular}{| c | cc  | cc| cc| cc| cc| }

\hline
Mesh  &{$L^1$ error} & Order  &  {$L^1$ error} & Order &  {$L^1$ error} & Order  &  {$L^1$ error} & Order \\
 \hline
  &  \multicolumn{2}{c|}{$P^1$ ELDG} & \multicolumn{2}{c|}{$P^1$ ELDG1}& \multicolumn{2}{c|}{$P^1$ ELDG2}& \multicolumn{2}{c|}{$P^1$ ELDG3}
  \\ \hline
    20 &     2.54E-03 &    --     &    2.42E-03 &    --    &    2.55E-03 &    --   &    2.49E-03 &    --   \\
    40 &     6.18E-04 &    2.03    &    5.97E-04 &    2.02    &    6.18E-04 &   2.04   &    5.99E-04 &   2.06   \\
    80 &     1.58E-04 &     1.96  &    1.55E-04 &    1.94  &    1.58E-04 &    1.96 &    1.55E-04 &    1.95   \\
   160 &     3.66E-05 &     2.11  &    3.62E-05 &    2.10  &    3.66E-05 &   2.11  &    3.62E-05 &    2.10   \\
\hline
  &  \multicolumn{2}{c|}{$P^2$ ELDG} & \multicolumn{2}{c|}{$P^2$ ELDG1}& \multicolumn{2}{c|}{$P^2$ ELDG2}& \multicolumn{2}{c|}{$P^2$ ELDG3}
  \\ \hline
    20 &    5.92E-05 & --        &    6.91E-05  &    --    &    6.01E-05  &    --   &    7.02E-05  &    --    \\
    40 &    7.48E-06 & --        &    7.83E-06 &    3.14    &    7.49E-06  &    3.00   &    7.81E-06  &    3.17   \\
    80 &    9.17E-07 &     3.03  &    9.29E-07 &   3.08  &    9.17E-07 &   3.03 &    9.29E-07  &    3.07   \\
   160 &    1.17E-07 &     2.97  &    1.18E-07 &   2.98   &    1.17E-07 &    2.97  &    1.18E-07  &   2.98   \\
\hline
\end{tabular}
\label{linear1d_spatial system for 'background scheme'}
\end{table}

\begin{table}[!ht]\footnotesize
\caption{1D wave equation with constant coefficient. $ u_{tt}  =  u_{xx} $ with initial condition $u(x,0) = \sin(x)$. $T=1$.
We use $CFL=0.3$ and $CFL=0.18$ with RK4 time discretization and post-processed technique for all $P^1$ and $P^2$ respectively. The error with post-processed technique for only $u^1=u_t$ was shown in this table.
  }
\centering
\begin{tabular}{| c | cc  | cc| cc| cc| cc| }

\hline
Mesh  &{$L^1$ error} & Order  &  {$L^1$ error} & Order &  {$L^1$ error} & Order  &  {$L^1$ error} & Order \\
 \hline
  &  \multicolumn{2}{c|}{$P^1$ ELDG} & \multicolumn{2}{c|}{$P^1$ ELDG1}& \multicolumn{2}{c|}{$P^1$ ELDG2}& \multicolumn{2}{c|}{$P^1$ ELDG3}
  \\ \hline
    20 &     2.26E-04 &    --     &    2.49E-04 &    --    &    2.38E-04 &    --   &    2.39E-04 &    --   \\
    40 &     2.36E-05 &    3.26    &    2.40E-05 &    3.38   &    2.40E-05 &   3.31   &    2.35E-05 &   3.35   \\
    80 &     2.66E-06 &     3.15  &    2.64E-06 &    3.18  &    2.67E-06 &    3.16 &    2.62E-06 &    3.16   \\
   160 &     3.15E-07 &     3.08  &    3.11E-07 &    3.08  &    3.15E-07 &   3.09  &    3.11E-07 &    3.08   \\
\hline
  &  \multicolumn{2}{c|}{$P^2$ ELDG} & \multicolumn{2}{c|}{$P^2$ ELDG1}& \multicolumn{2}{c|}{$P^2$ ELDG2}& \multicolumn{2}{c|}{$P^2$ ELDG3}
  \\ \hline
    20 &    2.15E-06 & --        &    2.27E-06 &    --    &    2.19E-06  &    --   &    2.28E-06  &    --    \\
    40 &    3.63E-08 & 5.89      &    3.86E-08 &   5.87    &    3.67E-08  &    5.90   &    3.89E-08  &  5.87    \\
    80 &    6.40E-10 &     5.83  &    6.79E-10 &   5.83  &    6.46E-10 &   5.83 &    6.84E-10  &    5.83   \\
   160 &    1.27E-11 &     5.66  &    1.33E-11 &   5.68   &    1.28E-11 &    5.66  &    1.34E-11  &   5.68   \\
\hline
\end{tabular}
\label{linear1d_spatial system for 'background scheme'post}
\end{table}


\begin{figure}[!ht]
\centering
\includegraphics[height=50mm]{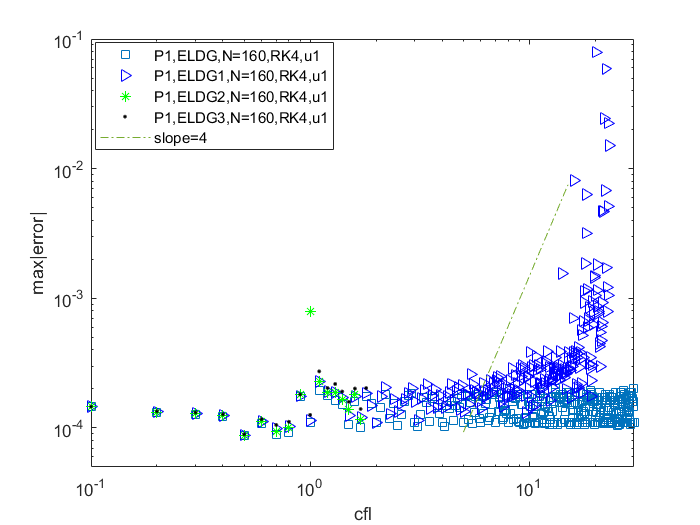}
\includegraphics[height=50mm]{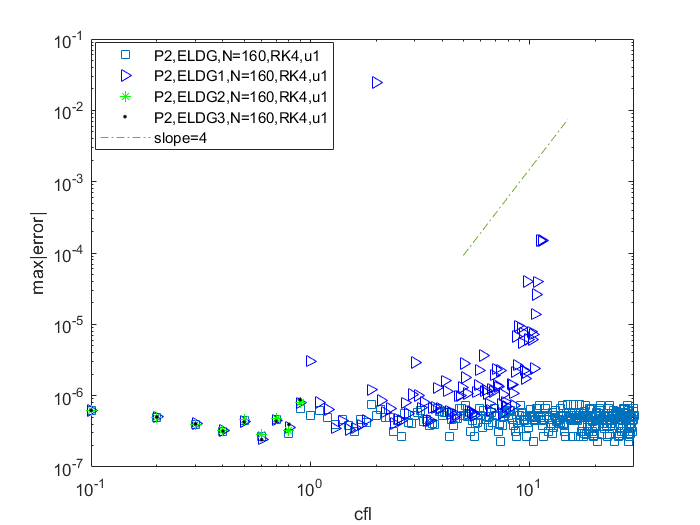}
\caption{The $L^\infty$ error versus $CFL$ of ELDG methods, ELDG1, ELDG2 and ELDG3 methods for 1D wave equation with constant coefficient:  $u_{tt}=u_{xx}$ with initial condition $u(x,0) = \sin(x)$. A long time simulation is performed with $T=100$ and mesh size $N=160$.  
}
\label{linear_stability1 for system}
\end{figure}

\begin{figure}[!ht]
\centering
\includegraphics[height=50mm]{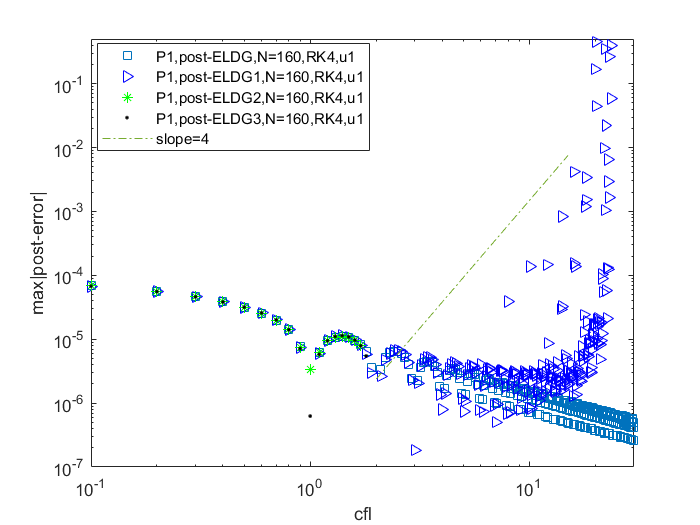}
\includegraphics[height=50mm]{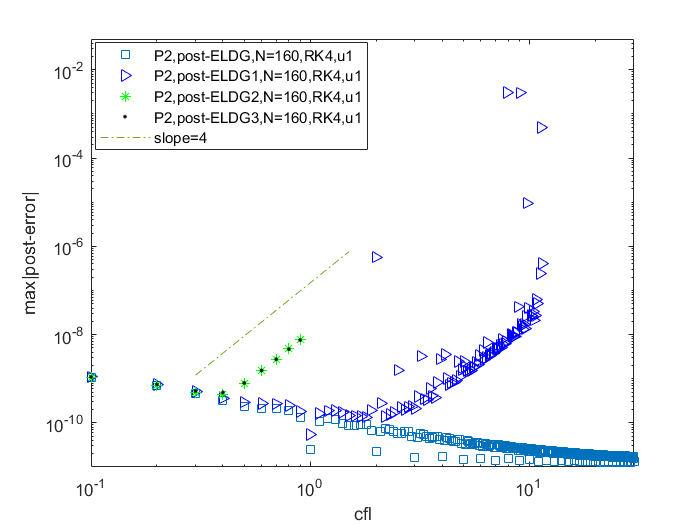}
\caption{The $L^\infty$ error versus $CFL$ of ELDG methods, ELDG1, ELDG2 and ELDG3 methods with post-processed technique for 1D wave equation with constant coefficient:  $u_{tt}=u_{xx}$ with initial condition $u(x,0) = \sin(x)$. A long time simulation is performed with $T=100$ and mesh size $N=160$.  
}
\label{linear_stability1 for systempost}
\end{figure}
\begin{figure}[!ht]
\centering
\includegraphics[height=50mm]{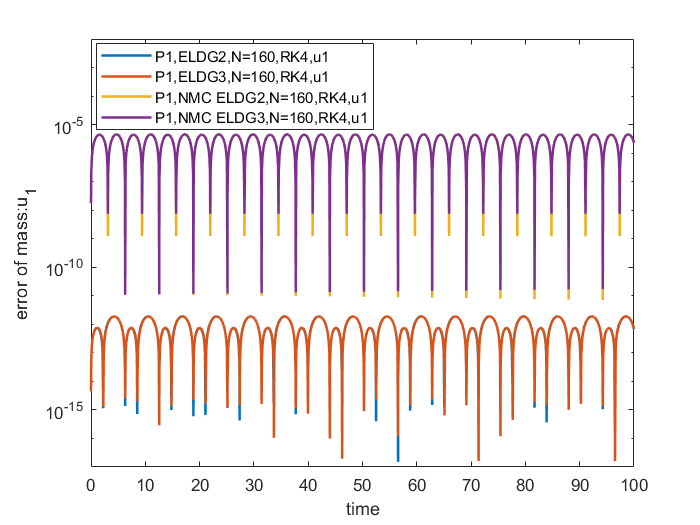}
\includegraphics[height=50mm]{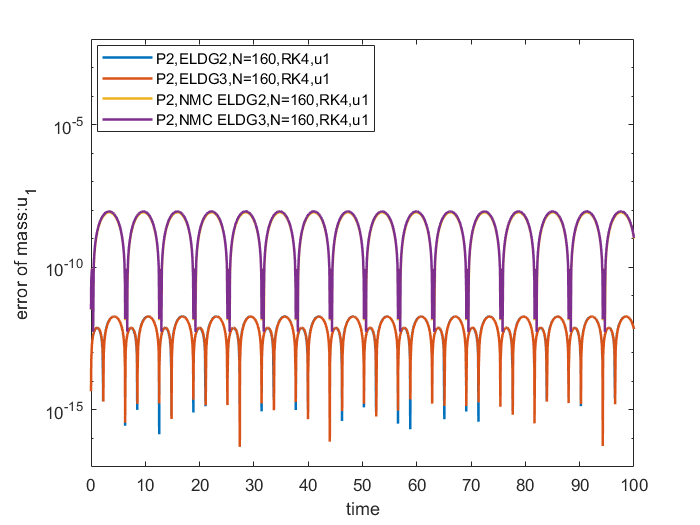}
\caption{The error of mass versus time of ELDG2, ELDG3, NMC ELDG2 and NMC ELDG3 methods for 1D wave equation with constant coefficient:  $u_{tt}=u_{xx}$ with initial condition $u(x,0) = \sin(x)$. A long time simulation $T=100$ is performed with meshes $N=160$, $CFL=0.1$ and RK4 time discretization.  
}
\label{linearmassconservesmoothinitial}
\end{figure}
\end{exa}

\begin{exa}
Then, we test the ELDG schemes for eq.~\eqref{2d_hcl} with a Gaussian initial condition $$u^1=exp(-\frac{x^2}{0.005}), u^2=0.$$ The computational domain is $[-1,1]$ with the periodic boundary conditions. The exact solutions $u^1=0.5(exp(-\frac{(x+t)^2}{0.005})+exp(-\frac{(x-t)^2}{0.005}))$ and $u^2=0.5(exp(-\frac{(x+t)^2}{0.005})-exp(-\frac{(x-t)^2}{0.005}))$ are the superposition of two Gaussian functions
 with a periodic extensions.
We plot $u_1$ ELDG3 with $P^1$ and $P^2$ numerical solutions at time $T=50.5$ for this system in Figure \ref{utt=uxxloglogRK4gaussiansolution320}.
We can observed that there is no significant phase difference with a long time simulation, meanwhile the dissipation can be improved by the mesh refinement and higher order spatial approximation.

\begin{figure}[!ht]
\centering
\includegraphics[height=50mm]{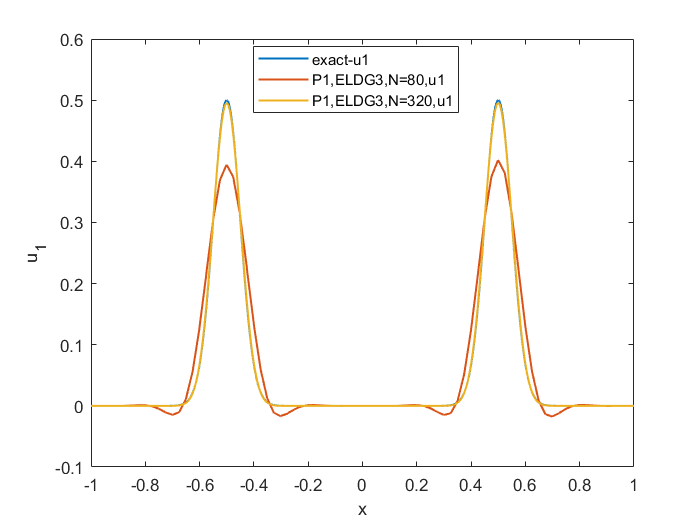}
\includegraphics[height=50mm]{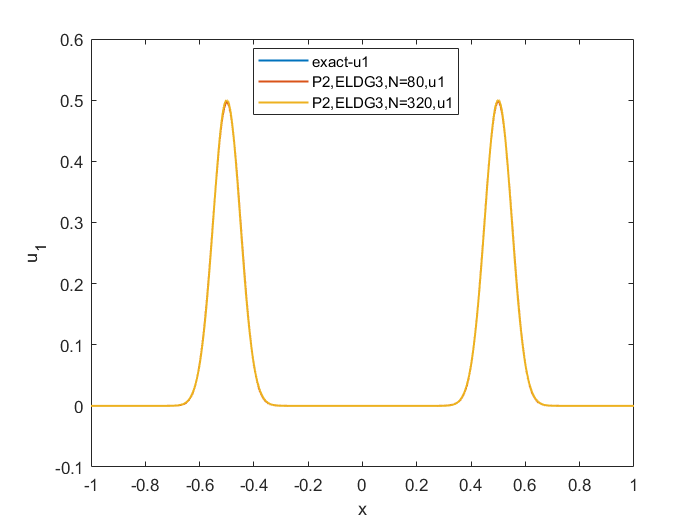}
\caption{Plots of the exact and numerical solutions $u_1$ at time $T=50.5$ of ELDG3 scheme for solving $u_{tt}=u_{xx}$ with Gaussian function initial condition. The mesh size of $N=80$ and $N=320$ are used. Left: $k=1$ ELDG3 with $CFL = 1.5$. Right: $k=2$ ELDG3 with $CFL = 0.9$.
}
\label{utt=uxxloglogRK4gaussiansolution320}
\end{figure}
\end{exa}

\begin{exa}
Next, we test the ELDG schemes for eq.~\eqref{2d_hcl}
on $[0, 2\pi]$ with the periodic boundary conditions and the following discontinuous initial condition
\begin{equation}\label{stepfunctioninitial}
\begin{aligned}
  &u^1_0(x)=\left\{
  \begin{aligned}
    &1., &if \ 0.95\pi \leq x\leq 1.05\pi,\\
    &0.5, \ &otherwise,
  \end{aligned}
  \right.\\
  &u^2_0(x)=1.
\end{aligned}
\end{equation}
The exact solutions $u^1$ and $u^2$ are discontinuous piecewise constants with moving discontinuities. It is a challenging test for controlling oscillations around discontinuities. We adopt a simple TVD limiter on background mesh at each RK stages with M = 0 in \cite{cockburn1989tvb} for all schemes.
As shown in Figure \ref{linear_stability1 for system}, the $CFL$ constraint with stability is slightly less than 1 for ELDG3 scheme. 
We plot the numerical solutions $u_1$ of ELDG3 scheme with $P^1$ and $P^2$, $CFL=0.9$ in Figure \ref{utt=uxxloglogRK4reimannsolution}.
It is found that oscillations are well controlled with the TVD limiter and ELDG method performs well for large time stepping size.
Moreover, we track the mass conservation of ELDG methods VS NMC ELDG methods for eq.~\eqref{2d_hcl} and present results in Figure \ref{reimannlinearmassconservestepinitialcfl}. It shows that the ELDG schemes maintain the mass conservation at the level of machine error, while the NMC ELDG schemes do not.

\begin{figure}[!ht]
\centering
\includegraphics[height=50mm]{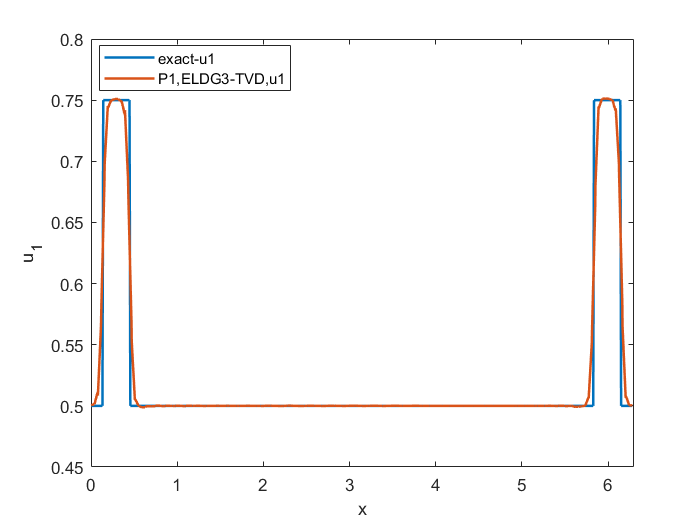}
\includegraphics[height=50mm]{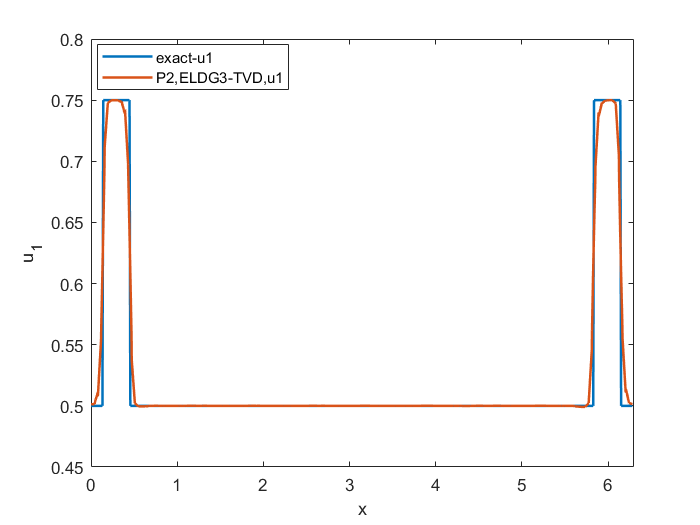}
\caption{Plots of the numerical solutions $u_1$ of ELDG3 scheme with TVD limiter for solving $u_{tt}=u_{xx}$ with step function initial condition. The final integration time T is 2.85. The mesh of $160$ is used. Left: $k=1$ ELDG3+TVDlimiter with $CFL = 1.5$. Right: $k=2$ ELDG3+TVDlimiter with $CFL = 0.9$.
}
\label{utt=uxxloglogRK4reimannsolution}
\end{figure}
%

\begin{figure}[!ht]
\centering
\includegraphics[height=50mm]{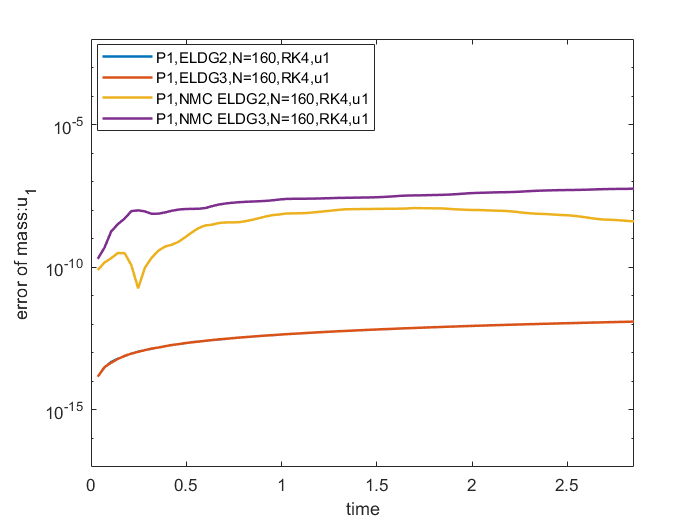}
\includegraphics[height=50mm]{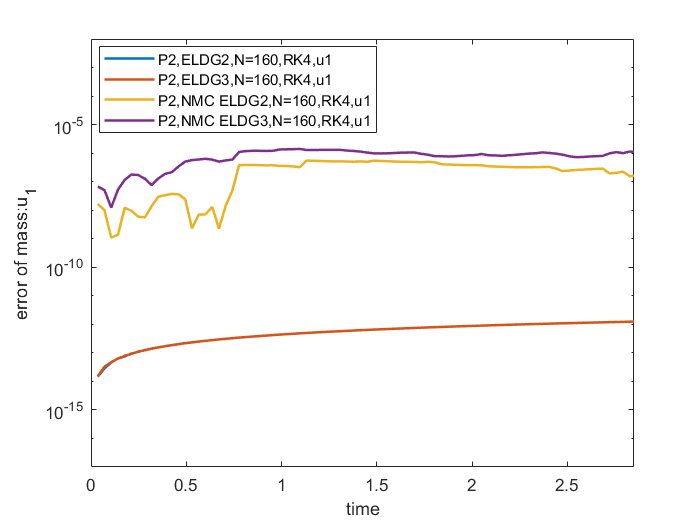}
\caption{The error of mass versus time of ELDG2, ELDG3, NMC ELDG2 and NMC ELDG3 methods with TVD limiter for 1D wave equation with constant coefficient:  $u_{tt}=u_{xx}$ with initial condition step function. $T=2.85$, $N=160$, $CFL=0.9$ and RK4 time discretization are performed for the simulation.  
}
\label{reimannlinearmassconservestepinitialcfl}
\end{figure}
\end{exa}

\begin{exa}
(1D wave equation with variable coefficient and source term.)
We consider the 1D wave equation eq.~\eqref{2d_hcl} with variable coefficient $a(x)=2+\sin(x)$ and exact solution $u(x,t)=\sin(x-2t)$ is periodic on $[0,2\pi]$. The source term is $f(x,t)=-4\sin(x-2t)+\sin(x-2t)(2+\sin(x))^2-2(2+\sin(x))\cos(x)\cos(x-2t)$. For computation, we choose mesh velocity $\nu^{(1)}_{j+\frac12}=a(x_{j+\frac12}), \nu^{(2)}_{j+\frac12}=-\nu^{(1)}_{j+\frac12}$ and exact eigenvectors with $a_p(x)=2+\sin(x)$ in \eqref{approxieigenvecterandarepresent}.

\begin{table}[!ht]\footnotesize
\caption{1D wave equation with variable coefficient and source term. $ u_{tt}  =  (a^2(x)u_x)_x+f(x,t) $ with initial condition $u(x,0) = \cos(x)$. $T=1$.
We use $CFL=0.1$ for $P^1$ and $P^2$ with RK4 time discretization. The error for only $u^1=u_t$ was shown in this table.
  }
\centering
\begin{tabular}{| c | cc  | cc| cc| cc| }

\hline
Mesh  &{$L^1$ error} & Order  &  {$L^2$ error} & Order &  {$L^\infty$ error} & Order \\
 \hline
  &  \multicolumn{6}{c|}{$P^1$ }
  \\ \hline
  20  &  5.20E-03 &  -- &  6.70E-03&
  -- &  2.39E-02 &  --\\
  40  &  1.32E-03 &  1.98 &  1.74E-03&
  1.94 &  6.54E-03 &  1.87\\
  80  &  3.28E-04 &  2.00 &  4.40E-04&
  1.99 &  1.68E-03 &  1.96\\
  160 &  8.11E-05 &  2.02 &  1.09E-04&
  2.01 &  4.16E-04 &  2.01\\

\hline
  &  \multicolumn{6}{c|}{$P^2$ }
  \\ \hline
 20 &  1.16E-04 &  --  & 1.61E-04&
   --  & 5.48E-04 &  -- \\
 40 &  1.48E-05 & 2.97 & 2.01E-05&
  3.00 & 6.70E-05 & 3.03\\
 80 &  1.88E-06 & 2.98 & 2.46E-06&
  3.03 & 7.84E-06 & 3.10\\
 160 &  2.31E-07 & 3.02 & 3.13E-07&
  2.98 & 1.04E-06 & 2.91\\

\hline

\end{tabular}
\label{1D general system for 'background scheme' and 'projection scheme'}
\end{table}

\begin{table}[!ht]\footnotesize
\caption{1D wave equation with variable coefficient and source term. $ u_{tt}  =  (a^2(x)u_x)_x+f(x,t) $ with initial condition $u(x,0) = \cos(x)$. $T=1$.
We use $CFL=0.1$ for $P^1$ and $P^2$ with post-processed technique and RK4. The error for only $u^1=u_t$ was shown in this table.
  }
\centering
\begin{tabular}{| c | cc  | cc| cc| cc| }

\hline
Mesh  &{$L^1$ error} & Order  &  {$L^2$ error} & Order &  {$L^\infty$ error} & Order \\
 \hline
  &  \multicolumn{6}{c|}{$P^1$ }
  \\ \hline
  20  &  9.10E-04 & -- & 1.07E-03&
  -- & 1.91E-03 & --\\
  40  &  1.07E-04 & 3.09 & 1.26E-04&
  3.09 & 2.26E-04 & 3.08\\
  80  &  1.29E-05 & 3.05 & 1.53E-05&
  3.05 & 2.75E-05 & 3.04\\
  160 &  1.58E-06 & 3.02 & 1.88E-06&
  3.02 & 3.39E-06 & 3.02\\


\hline
  &  \multicolumn{6}{c|}{$P^2$ }
  \\ \hline
  20  &  5.34E-06 & -- & 6.39E-06&
  -- & 1.56E-05 & --\\
  40  &  8.95E-08 & 5.90 & 1.03E-07&
  5.96 & 2.82E-07 & 5.79\\
  80  &   1.73E-09 & 5.69 & 1.94E-09&
  5.73 & 3.32E-09 & 6.41\\
  160 &  6.62E-11 & 4.71 & 7.59E-11&
  4.67 & 1.46E-10 & 4.51\\
\hline

\end{tabular}
\label{1D general system for 'background scheme' and 'projection scheme'post}
\end{table}

\begin{figure}[!ht]
\centering
\includegraphics[height=50mm]{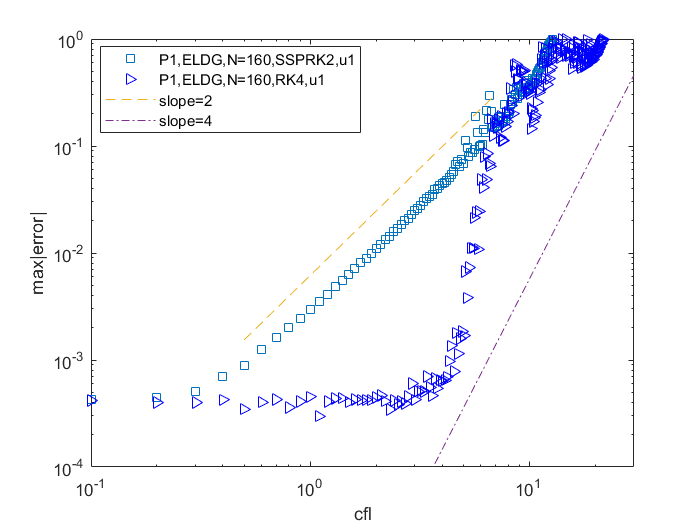}
\includegraphics[height=50mm]{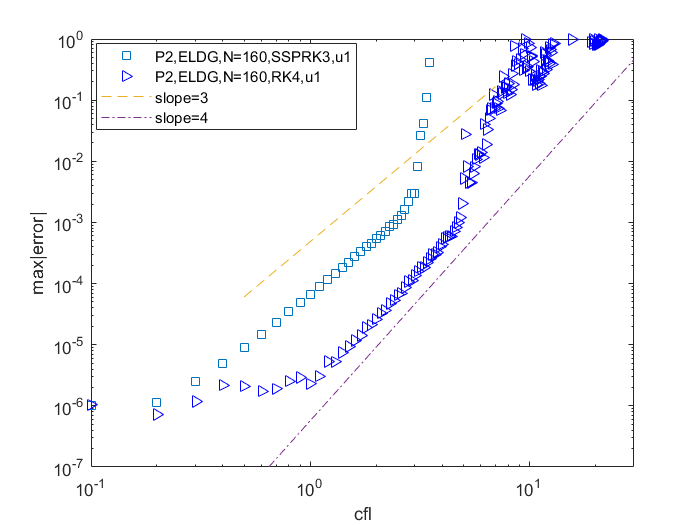}
\caption{The $L^\infty$ error versus $CFL$ of ELDG method for $ u_{tt}  =  (a^2(x)u_x)_x+f(x,t) $ with initial condition $u(x,0) = \cos(x)$. $T=1$.    $\Delta t=CFL\Delta x$.
}
\label{sin1d_time}
\end{figure}

\begin{figure}[!ht]
\centering
\includegraphics[height=50mm]{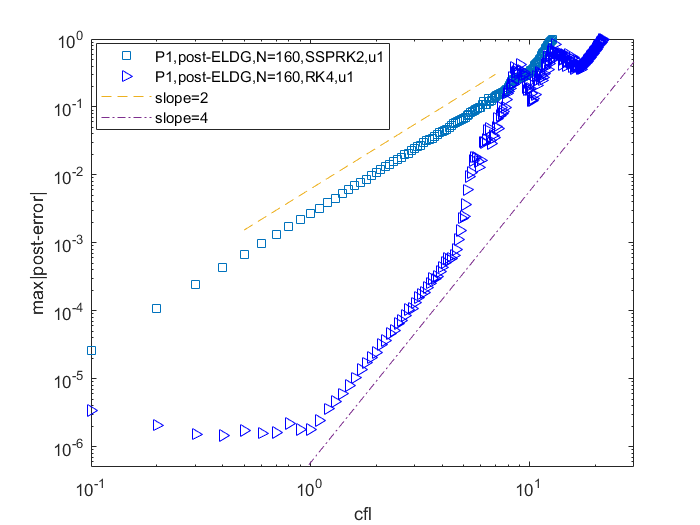}
\includegraphics[height=50mm]{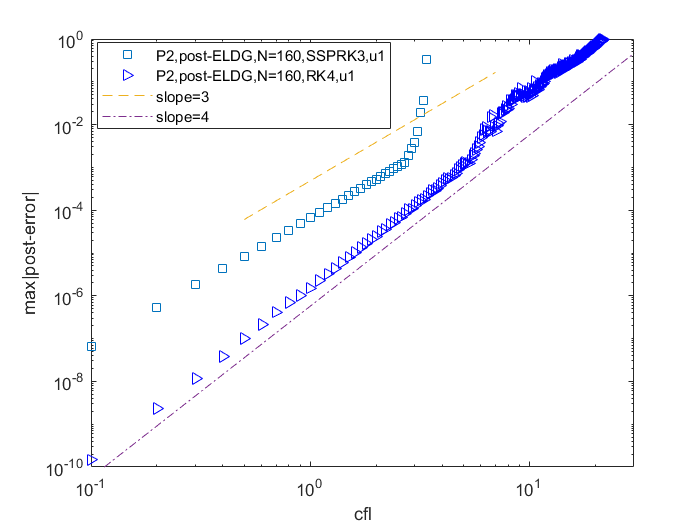}
\caption{The $L^\infty$ error versus $CFL$ of ELDG method with post-processed technique for $ u_{tt}  =  (a^2(x)u_x)_x+f(x,t) $ with initial condition $u(x,0) = \cos(x)$. $T=1$.    $\Delta t=CFL\Delta x$.
}
\label{sin1d_timepost}
\end{figure}


The expected optimal spatial accuracies of the ELDG methods without and with post-processing technique are shown in Table \ref{1D general system for 'background scheme' and 'projection scheme'} and Table \ref{1D general system for 'background scheme' and 'projection scheme'post} respectively.
In Figure \ref{sin1d_time} and \ref{sin1d_timepost}, we plot the $L^\infty$ error versus $CFL$ of ELDG methods without and with post-processing technique respectively. The following observations are made:
(1) The high order accuracy of the RK method reduce the error magnitude when large time stepping size is used; (2) The ELDG methods with RK4 time discretization perform well around and before $CFL=1$, which is well above the stability constraint of the RK DG method $1/(2k+1)$ for $P^k$ approximations. (3) After $CFL=1$ and before stability constraint of the method, the temporal convergence order is observed to be consistent with the order of RK discretization; (4) The ELDG methods with post-processing technique have smaller error magnitute than those without post-processing.

\end{exa}

\begin{exa}
(Two-dimensional linear system with constant coefficient matrices.)
The second order wave equation $u_{tt}=u_{xx}+u_{yy}$ can be written as the following first order linear hyperbolic system:
\begin{equation}
\begin{cases}
\left( \begin{matrix}
  u\\
  v
\end{matrix}\right)_t+\left( \begin{matrix}
  -1 &0\\
  0 &1
\end{matrix}\right)\left( \begin{matrix}
  u\\
  v
\end{matrix}\right)_x+\left( \begin{matrix}
  0 &-1\\
  -1 &0
\end{matrix}\right)\left( \begin{matrix}
  u\\
  v
\end{matrix}\right)_y=\left( \begin{matrix}
  0\\
  0
\end{matrix}\right),\\
u(x,y,0) = \frac{1}{2\sqrt{2}}\sin(x+y)-\frac{1}{2\sqrt{2}}\cos(x+y),\\
v(x,y,0) = \frac{\sqrt{2}-1}{2\sqrt{2}}\sin(x+y)+\frac{\sqrt{2}+1}{2\sqrt{2}}\cos(x+y)
\end{cases}
\label{eqn:2d_wave}
\end{equation}
with period boundary conditions in both $x$ and $y$ directions. The exact solution is
\begin{equation}
\begin{cases}
u(x,y,t) = \frac{1}{2\sqrt{2}}\sin(x+y+\sqrt{2}t)-\frac{1}{2\sqrt{2}}\cos(x+y-\sqrt{2}t),\\
v(x,y,t) = \frac{\sqrt{2}-1}{2\sqrt{2}}\sin(x+y+\sqrt{2}t)+\frac{\sqrt{2}+1}{2\sqrt{2}}\cos(x+y-\sqrt{2}t).
\end{cases}
\label{eqn:2d_wave_exactsolution}
\end{equation}
We notice that the two matrices in equation~\eqref{eqn:2d_wave} don't commute, thus the linear system can not be reduced to 2-D scalar problems.
We test accuracy for $Q^k$ ELDG methods with RK4 and 4th splitting method \cite{cai2016conservative} at $T=1$ for $k = 1, 2$ with $CFL = 0.1 $ in Table \ref{2Dconstantcoefsystemfor'background scheme'}. As expected, the $(k + 1)$th order convergence is observed for these methods. We plot the $L^\infty$ error versus $CFL$ of ELDG methods with $Q^1$ (left) and $Q^2$ (right) polynomial spaces for this case with Strang splitting and 4th order splitting in Figure \ref{constant2D}, which shows that second and forth order splitting errors are dominant when time-stepping sizes are large enough. The CFL constraint with stability for ELDG method is larger than that for general RK DG method when high order time discretization is applied.

\begin{table}[!ht]\footnotesize
\caption{Two-dimensional linear system with constant coefficient matrices. $Q^k$ EL DG methods (k = 1, 2) with RK4 and 4th splitting time discretization methods for \eqref{eqn:2d_wave} with the smooth initial condition at T = 1. $CFL = 0.1$.
  }
\centering
\begin{tabular}{| c | cc  | cc| cc| cc| }

\hline
Mesh  &{$L^1$ error} & Order  &  {$L^2$ error} & Order &  {$L^\infty$ error} & Order \\
 \hline
  &  \multicolumn{6}{c|}{$Q^1$}
  \\ \hline
  $20^2$  &  8.03E-04 & -- & 9.47E-04&
  -- & 1.85E-03 & --\\
  $40^2$  &  2.16E-04 & 1.89 & 2.50E-04&
  1.92 & 4.56E-04 & 2.02\\
  $80^2$  &  5.57E-05 & 1.96 & 6.40E-05&
  1.97 & 1.13E-04 & 2.02\\
  $160^2$ &  1.43E-05 & 1.96 & 1.64E-05&
  1.97 & 2.84E-05 & 1.99\\


\hline
  &  \multicolumn{6}{c|}{$Q^2$}
  \\ \hline
  $20^2$  &  1.70E-04 & -- & 1.90E-04&
  -- & 3.12E-04 & --\\
  $40^2$  &  2.21E-05 & 2.95 & 2.47E-05&
  2.94 & 4.14E-05 & 2.91\\
  $80^2$  &   2.75E-06 & 3.00 & 3.08E-06&
  3.00 & 5.21E-06 & 2.99\\
 $160^2$ &  3.38E-07 & 3.02 & 3.80E-07&
  3.02 & 6.45E-07 & 3.01\\
\hline

\end{tabular}
\label{2Dconstantcoefsystemfor'background scheme'}
\end{table}

\begin{figure}[!ht]
\centering
\includegraphics[height=50mm]{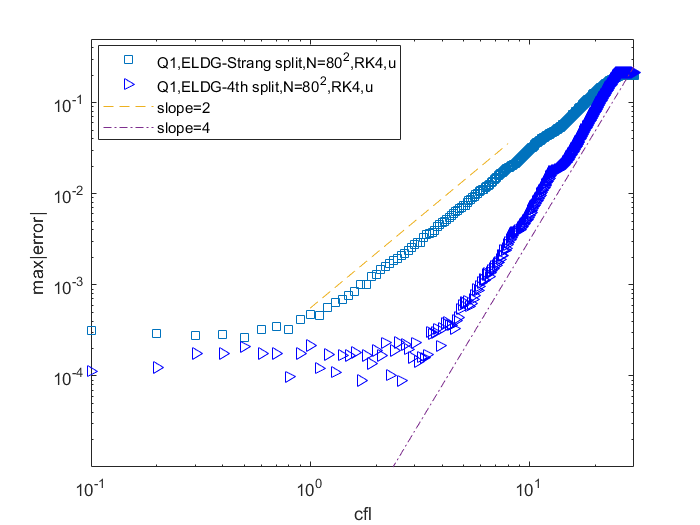}
\includegraphics[height=50mm]{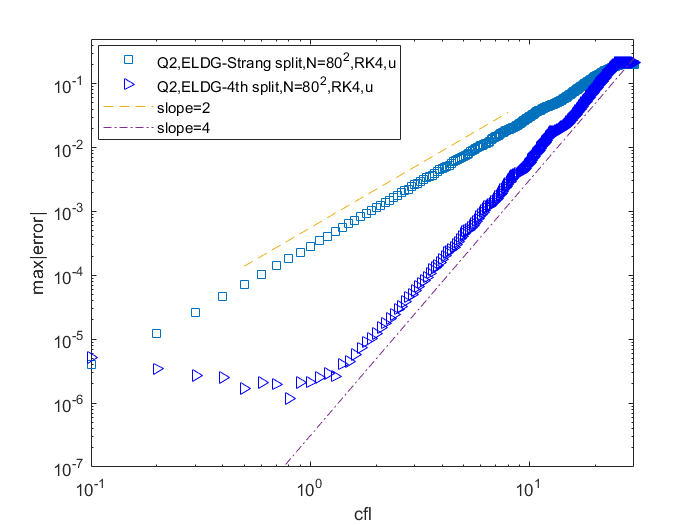}
\caption{The $L^\infty$ error versus $CFL$ of ELDG method with Strang splitting and 4th splitting, RK4 time discretization for \eqref{eqn:2d_wave}. $T=1$.
}
\label{constant2D}
\end{figure}
\end{exa}

\begin{exa}
(Two-dimensional linear system with variable coefficient matrices.)
The second order wave equation $u_{tt}=(a^2(x,y)u_x)_x+(b^2(x,y)u_y)_y$ can be written as the following first order linear hyperbolic system:
\begin{equation}
\begin{cases}
\left( \begin{matrix}
  u_1\\
  u_2\\
  u_3
\end{matrix}\right)_t+\left[\left( \begin{matrix}
  0 &-(a(x,y))^2 &0\\
  -1 &0 &0\\
  0 &0 &0
\end{matrix}\right)\left( \begin{matrix}
  u_1\\
  u_2\\
  u_3
\end{matrix}\right)\right]_x+\left[\left( \begin{matrix}
  0 &0 &-(b(x,y))^2)\\
  0 &0 &0\\
  -1 &0 &0
\end{matrix}\right)\left( \begin{matrix}
  u_1\\
  u_2\\
  u_3
\end{matrix}\right)\right]_y=\left( \begin{matrix}
  0\\
  0\\
  0
\end{matrix}\right).
\end{cases}
\label{eqn:2d_wavevariable}
\end{equation}
The simplest case is to choose $a(x,y)=b(x,y)=1$, computational domain $[-1,1]\times[-1,1]$ with periodic boundary condition and the Gaussian function initial condition:
\begin{equation}
\begin{cases}
u_1(x,y,0) = exp(-\frac{x^2+y^2}{0.005}),\\
u_2(x,y,0) = 0,\\
u_3(x,y,0) = 0.
\end{cases}
\end{equation}
For this example, we show the numerical ELDG $Q^2$ solution $u_1$ at times $T=0.5; 1; 1.5; 2$ in Figure \ref{solutiongaussianTcontour}.

\begin{figure}[!ht]
\centering
\includegraphics[height=50mm]{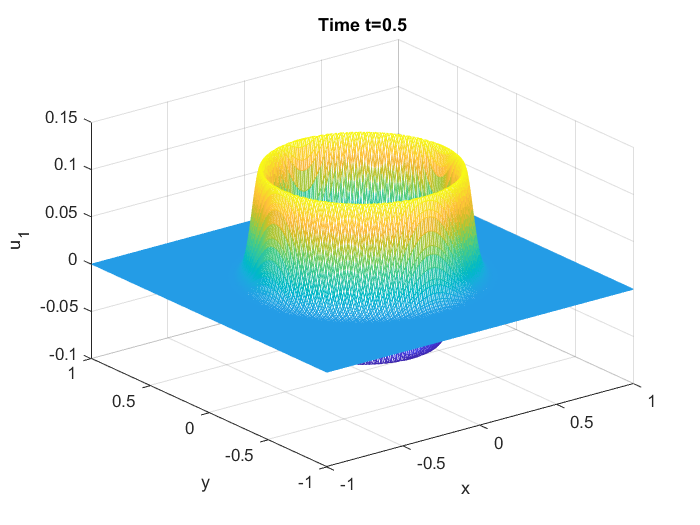}
\includegraphics[height=50mm]{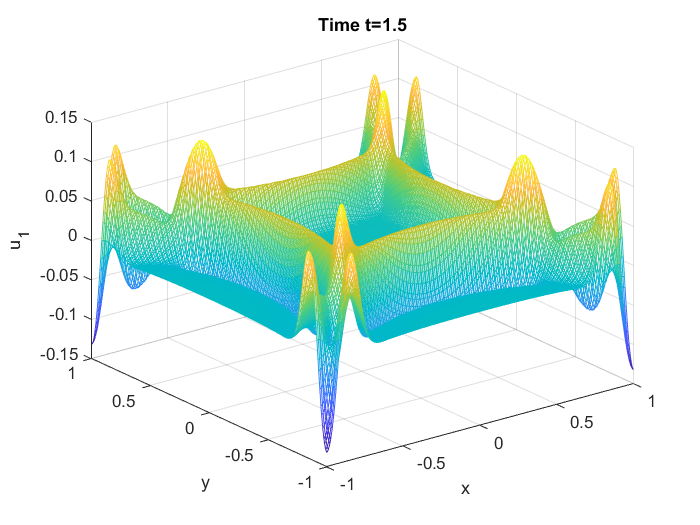}
\includegraphics[height=50mm]{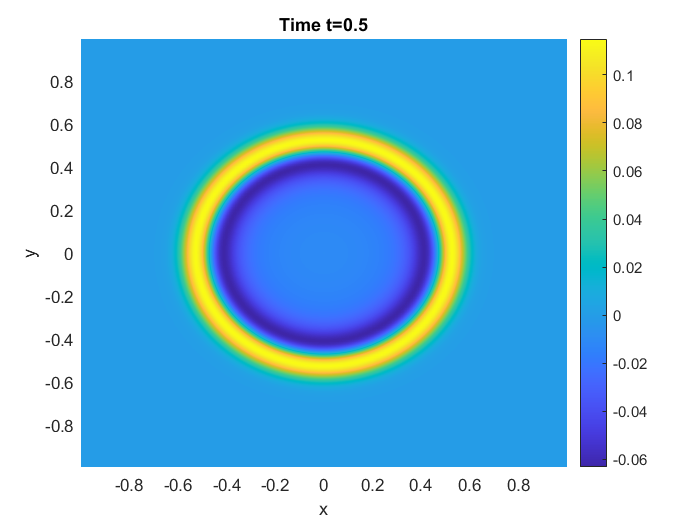}
\includegraphics[height=50mm]{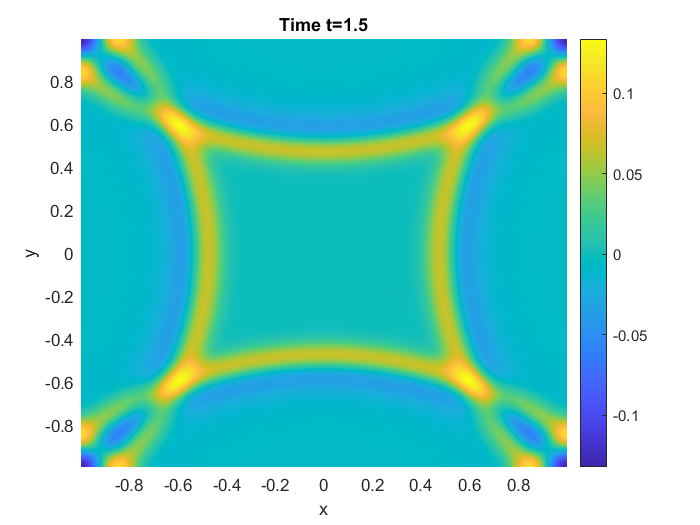}
\includegraphics[height=50mm]{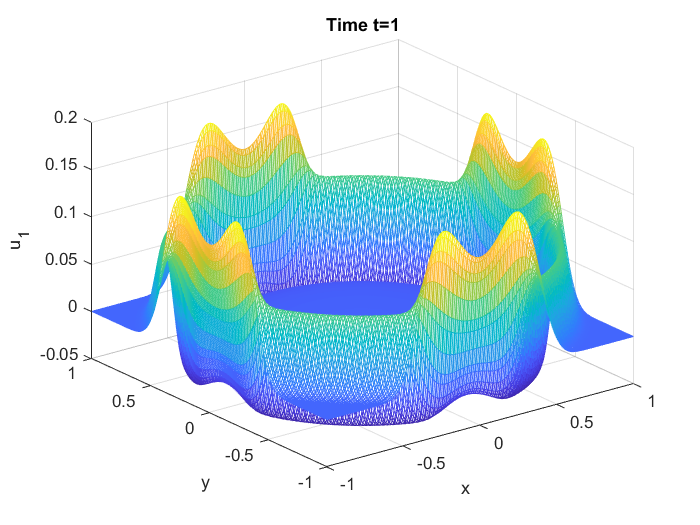}
\includegraphics[height=50mm]{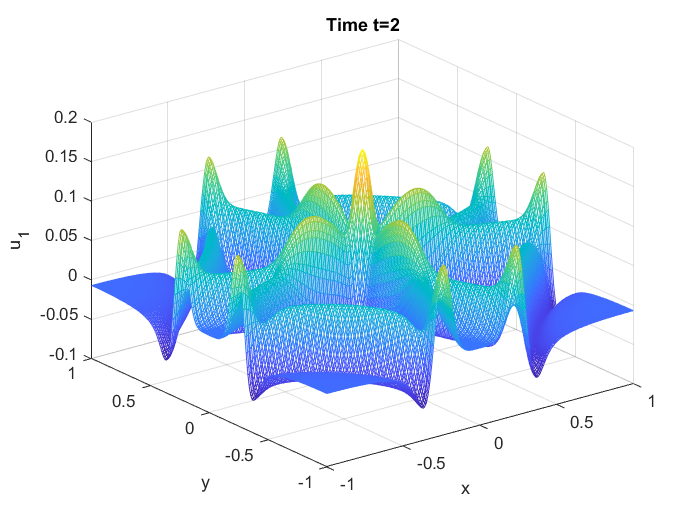}
\includegraphics[height=50mm]{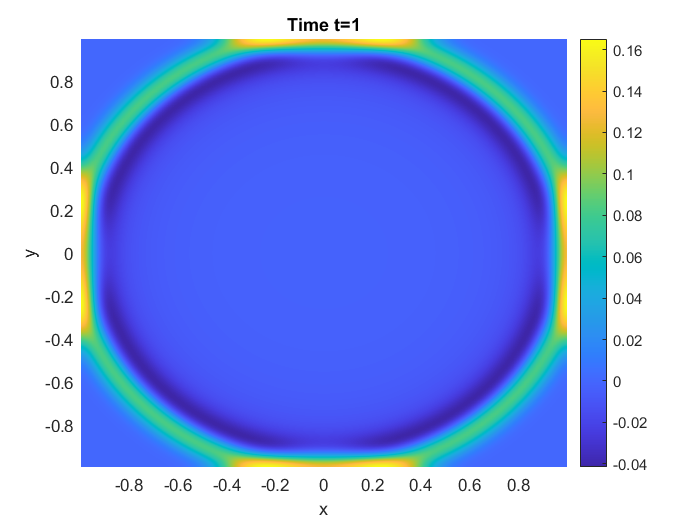}
\includegraphics[height=50mm]{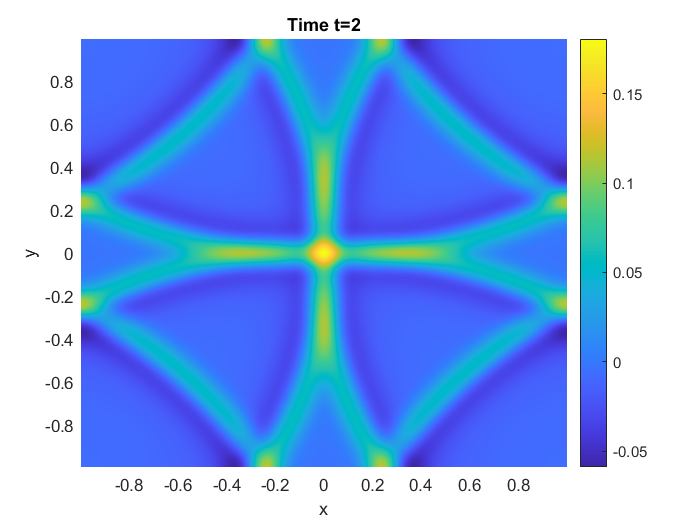}
\caption{Plots of the ELDG numerical solutions $u_1$ and their contour plots at $T=0.5; 1; 1.5; 2$ for 2D wave equation $u_{tt}=u_{xx}+u_{yy}$ with Gaussian function initial condition. The mesh of $80\times80$ is used with 4th splitting method and RK4 time discretization, $CFL=1$.
}
\label{solutiongaussianTcontour}
\end{figure}

Next we consider the system \eqref{eqn:2d_wavevariable} with the initial condition
\begin{equation}
\begin{cases}
u_1(x,y,0) = 2\cos(x+y),\\
u_2(x,y,0) = \cos(x+y),\\
u_3(x,y,0) = \cos(x+y),
\end{cases}
\end{equation}
where $a(x,y)=1+0.5\sin(x+y)$, $b(x,y)=\sqrt{(4-(1+0.5\sin(x+y))^2)}$ and the boundary condition is periodic in both $x$ and $y$ directions. The exact solution is
\begin{equation}
\begin{cases}
u_1(x,y,t) = 2\cos(x+y+2t),\\
u_2(x,y,t) = \cos(x+y+2t),\\
u_3(x,y,t) = \cos(x+y+2t).
\end{cases}
\label{eqn:2d_wave_exactsolutionvariable}
\end{equation}
We report the spatial accuracy of $Q^k$ ELDG methods in Table \ref{2Dvariablecoefsystemfor'background scheme'}. The expected optimal convergence is observed. We plot the $L^\infty$ error versus $CFL$ of ELDG methods in Figure \ref{variable2D}. The ELDG methods perform as well as that for the linear system with constant coefficient matrices, and the CFL allowed with stability is much larger than that of the RK DG method.
\begin{table}[!ht]\footnotesize
\caption{Two-dimensional linear system with variable matrices. $Q^k$ ELDG methods (k = 1, 2) with RK4 and 4th splitting time discretization methods for \eqref{eqn:2d_wavevariable} with the smooth initial condition at T = 0.1, $CFL = 0.1$.
  }
\centering
\begin{tabular}{| c | cc  | cc| cc| cc| }

\hline
Mesh  &{$L^1$ error} & Order  &  {$L^2$ error} & Order &  {$L^\infty$ error} & Order \\
 \hline
  &  \multicolumn{6}{c|}{$Q^1$}
  \\ \hline
  $20^2$  &  1.89E-03 & -- & 2.46E-03&
  -- & 5.29E-03 & --\\
  $40^2$  &  4.58E-04 & 2.05 & 6.01E-04&
  2.03 & 1.26E-03 & 2.07\\
  $80^2$  &  1.13E-04 & 2.02 & 1.50E-04&
  2.00 & 3.10E-04 & 2.02\\
  $160^2$ &  2.81E-05 & 2.01 & 3.73E-05&
  2.01 & 7.67E-05 & 2.01\\


\hline
  &  \multicolumn{6}{c|}{$Q^2$}
  \\ \hline
  $20^2$  &  2.95E-04 & -- & 3.66E-04&
  -- & 8.43E-04 & --\\
  $40^2$  &  4.06E-05 & 2.86 &  4.91E-05&
  2.90 & 1.01E-04 & 3.07\\
  $80^2$  &   5.15E-06 & 2.98 & 6.20E-06&
  2.99 & 1.30E-05 & 2.95\\
 $160^2$ &  6.49E-07 & 2.99 & 7.80E-07&
  2.99 & 1.63E-06 & 2.99\\
\hline

\end{tabular}
\label{2Dvariablecoefsystemfor'background scheme'}
\end{table}

\begin{figure}[!ht]
\centering
\includegraphics[height=50mm]{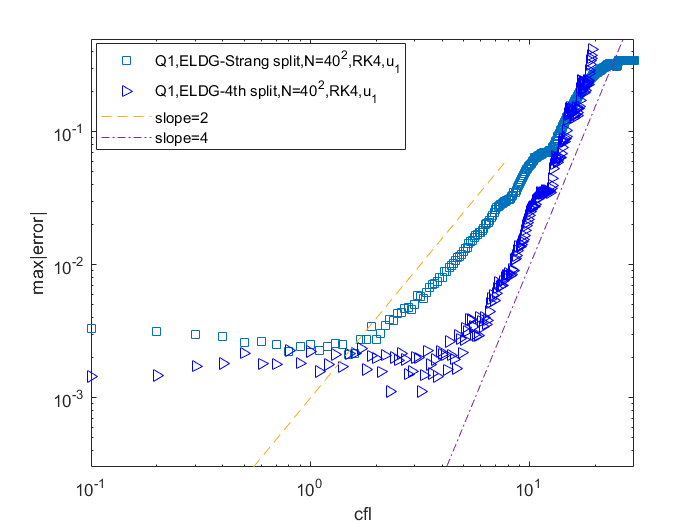}
\includegraphics[height=50mm]{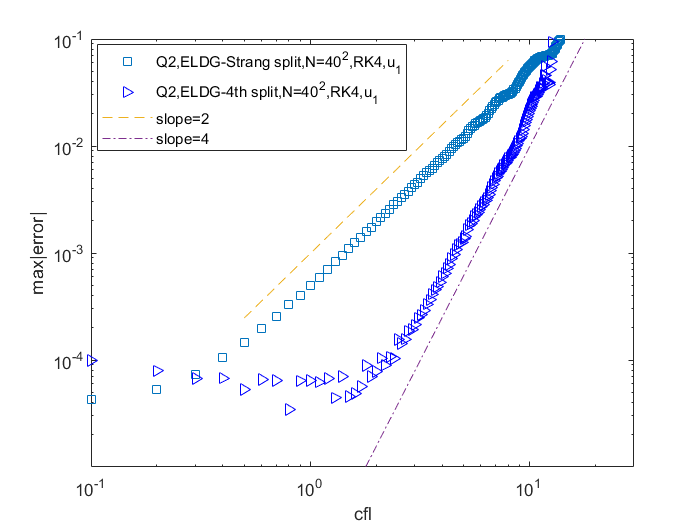}
\caption{The $L^\infty$ error versus $CFL$ of ELDG method with Strang splitting and 4th splitting, RK4 time discretization for \eqref{eqn:2d_wavevariable}. $T=1$, mesh size$40^2$.
}
\label{variable2D}
\end{figure}

\end{exa}

\section{Conclusion} \label{section:conclusion}

In this paper, we have developed a mass-conservative Eulerian-Lagrangian discontinuous Galerkin (ELDG) method for wave equations written in the form of hyperbolic systems.
The new framework track the information of each characteristic family by the corresponding characteristic region, and combine in a mass-conservative fashion; The method inherits advantages in stability under large time stepping sizes, and in mass conservation, compactness and high order accuracy.
These advantages are numerically verified by extensive numerical tests for 1D and 2D linear wave equations.
Future works include further theoretic development and application of ELDG methods for nonlinear hyperbolic problems.

\section{Appendix (The notation of test function)}\label{section:appendix}
We give some notations which is used in the implementation of the fully discrete ELDG scheme with RK time discretization.
For convenience, we only give the definitions related to $\Omega_j^{(1)}$ below because we can similarly get the definitions related to $\Omega_j^{(2)}$. As scalar case, we take test function $\psi^{(1)}(x,t)$ as $\psi^{(1)}_{j,m}(x,t)=\Psi_{j,m}(x-\alpha^{(1)} (t-t^{n+1}))$ in the adjoint problem which is a set of basis of $P^k(\tilde{I}^{(1)}_j(t))$, where $\tilde{I}^{(1)}_j(t)$ is donated by a domain related $\lambda^{(1)}$ as $\tilde{I}_j(t)$ in scalar problem.
Here, we also take $\Psi_{j,m}$ as orthogonal basis on $I_j$, and let
   \begin{align}\label{representation of basis}
    &u_h^{1,(1)}(x,t)=\sum_{l=0}^k \hat{u}^{1,(1);(l)}_j(t)\psi^{(1)}_{j,l}(x,t), \quad \mbox{on} \quad \tilde{I}^{(1)}_j(t),
   \end{align}
   where $ \hat{u}^{1,(1);(l)}$ are coefficients for the basis. Let $\hat{U}^{1,(1)}_j(t) = (\hat{u}_j^{1,(1);(0)}(t), \cdots, \hat{u}_j^{1,(1);(k)}(t))^T$ be the coefficient vector of size $(k+1)\times 1$. Then we have $$[\int_{\tilde{I}^{(1)}_j(t)}u_h^{1,(1)}(x,t)\psi^{(1)}_{j,0}(x,t)dx,...,\int_{\tilde{I}^{(1)}_j(t)}u_h^{1,(1)}(x,t)\psi^{(1)}_{j,k}(x,t)dx]^T=\hat{U}^{1,(1)}_j(t), \forall t\in [t^n,t^{n+1}].$$ $\hat{U}^{2,(1)}_j(t)$ can be similarly defined.
Similar definition can be made to $\Omega^{(2)}_j$ and $\psi^{(2)}(x,t)$ for the second characteristics family.

\bibliographystyle{abbrv}
\bibliography{refer17}

\end{document}